\documentclass[11pt,twoside]{article}

\usepackage{epsfig,amsfonts,color}
\usepackage{amsmath, bbm, mathabx}
\usepackage[normalem]{ulem}

\usepackage{amssymb, palatino, geometry,url}
\usepackage{algorithmic}
\usepackage[noresetcount,lined,boxed]{algorithm2e} % ... for algorithms

\usepackage[colorlinks=true,linkcolor=blue,citecolor=blue,urlcolor=blue]{hyperref}
\usepackage{subcaption}
\usepackage{amsthm}
\usepackage[utf8]{inputenc}
\usepackage{textgreek}

\usepackage{fancyhdr}
\usepackage{lineno}
\usepackage{float}
\usepackage[section]{placeins}

\usepackage{listings}
\DeclareCaptionFont{white}{\color{white}}
\DeclareCaptionFormat{listing}{{\parbox{\dimexpr\textwidth-2\fboxsep\relax}{#1#2#3}}}
\usepackage[dvipsnames]{xcolor}

\lstdefinelanguage{julia}{
    alsoletter={0123456789!.:+-^?*~><='/},
    keywords={},
    keywordsprefix=:,
    keywords = [2]{function, end, begin, if, elseif, else, using, for},
    keywords = [3]{InfiniteModel, MvNormal, optimize!, objective_value, Infinite, zeros,
                   DomainRestrictions, set_value, enumerate, value, value., length, supports,
                   OrthogonalCollocation, sum, add_supports, Uniform},
    keywords = [4]{0, 1, 2, 5, 10, 100, 1000, 0.95, 0.8, 0.2, 101, 20, 1:2, 1:5, 1:10, 15, 
                   4, 0.09, 2.1, -10, true, false},
    keywords = [5]{@variable, @constraint, @objective, @parameter_function, 
                   @infinite_parameter, @finite_parameter, @variables, @constraints},
    keywords = [6]{+, -, ^, ?, *, >, <, =, <=, >=, ==, .<=, .>=, .==, .+, .-, =>, in, /},
    keywordsprefix=:,
    keywordstyle=\color{red},
    keywordstyle=[2]\color{violet},
    keywordstyle=[3]\color{blue},
    keywordstyle=[4]\color{orange},
    keywordstyle=[5]\color{purple},
    keywordstyle=[6]\color{teal},
    sensitive=true,
    morecomment=[l]{\#},
    morestring=[b]"
}

\title{A Unifying Modeling Abstraction for \\ Infinite-Dimensional Optimization}
\author{Joshua L. Pulsipher, Weiqi Zhang, Tyler J. Hongisto, and Victor M. Zavala\thanks{Corresponding Author: victor.zavala@wisc.edu}\\
	{\small Department of Chemical and Biological Engineering}\\
	{\small \;University of Wisconsin-Madison, 1415 Engineering Dr, Madison, WI 53706, USA}}
\date{}

\begin{document}
	
\maketitle

\begin{abstract}
Infinite-dimensional optimization (InfiniteOpt) problems involve modeling components (variables, objectives, and constraints) that are functions defined over infinite-dimensional domains. Examples include continuous-time dynamic optimization (time is an infinite domain and components are a function of time), PDE optimization problems (space and time are infinite domains and components are a function of space-time), as well as stochastic and semi-infinite optimization (random space is an infinite domain and components are a function of such random space). InfiniteOpt problems also arise from combinations of these problem classes (e.g., stochastic PDE optimization). Given the infinite-dimensional nature of objectives and constraints, one often needs to define appropriate quantities (measures) to properly pose the problem. Moreover, InfiniteOpt problems often need to be transformed into a finite dimensional representation so that they can be solved numerically. In this work, we present a unifying abstraction that facilitates the modeling, analysis, and solution of InfiniteOpt problems. The proposed abstraction enables a general treatment of infinite-dimensional domains and provides a measure-centric paradigm to handle associated variables, objectives, and constraints. This abstraction allows us to transfer techniques across disciplines and with this identify new, interesting, and useful modeling paradigms (e.g., event constraints and risk measures defined over time domains). Our abstraction serves as the backbone of an intuitive {\tt Julia}-based modeling package that we call {\tt InfiniteOpt.jl}. We demonstrate the developments using diverse case studies arising in engineering. 
\end{abstract}


\section{Introduction} \label{sec:intro}

Infinite-dimensional optimization (InfiniteOpt) problems contain parameters that live in infinite-dimensional domains (e.g., time, space, random)  \cite{devolder2010solving}; the components of these problems (variables, objectives, and constraints) are parameterized over these domains and are thus infinite-dimensional functions (they form manifolds and surfaces).  A classical example of an InfiniteOpt problem is continuous-time dynamic optimization \cite{bertsekas1995dynamic}; here, the control trajectory is a function of time and time is a parameter that lives in an infinite-dimensional (continuous) domain. This formulation contrasts with that of a discrete-time dynamic optimization problem, in which the control trajectory is a collection of values defined over a finite set of times (domains are finite). Given the infinite-dimensional nature of variables, objectives, and constraints, one requires specialized techniques to define an InfiniteOpt problem properly. This is done by using {\em measures}, which are operators that summarize/collapse an infinite-dimensional object into a scalar quantity. For instance, in dynamic optimization,  one often minimizes the integral of the cost over the time domain and, in stochastic optimization, one often minimizes the expected value or variance of the cost. Measures are thus a key modeling element of InfiniteOpt problems that help manipulate the {\em shape} of infinite-dimensional objects to achieve desired outcomes (e.g., minimize peak/extreme costs or satisfy constraints with high probability). InfiniteOpt problems also often contain differential operators that dictate how components evolve over their corresponding domains; these operators often appear in problems that include differential and algebraic equations (DAEs) and partial differential equations (PDEs). 
\\

InfiniteOpt problems encompass a wide range of classical fields such as dynamic optimization \cite{biegler2007overview}, PDE optimization \cite{hinze2008optimization}, stochastic optimization \cite{birge2011introduction}, and semi-infinite optimization \cite{stein2003solving}. One also often encounters InfiniteOpt problems that are obtained by combining infinite-dimensional domains (e.g., space, time, and random domains). This situation arises, for instance, in stochastic dynamic optimization (e.g., stochastic optimal control) problems and in optimization problems with stochastic PDEs. In these problems, one needs to define measures that summarize modeling objects that are defined over the multiple domains (e.g., minimize the space-time integral of the cost or minimize the expected value of the time integral of the cost). InfiniteOpt problems appear in applications such as continuous-time model predictive control and moving horizon estimation  \cite{rawlings2000tutorial, qin2003survey}, design under uncertainty \cite{stankiewicz2000process, straub1993design}, portfolio planning \cite{ccakmak2006portfolio, dentcheva2006portfolio}, parameter estimation for differential equations \cite{shin2019scalable, biegler2003large}, reliability analysis \cite{pulsipher2020measuring, suksuwan2018optimization}, optimal power flow \cite{roald2015optimal, lan2018modeling}, and dynamic design of experiments \cite{georgakis2013design, asprey2002designing}.  
\\

The infinite-dimensional nature of modeling objects make InfiniteOpt problems challenging to solve \cite{nicholson2018pyomo, vanderbei2020linear, nocedal2006numerical}. Specifically,  these problems often need to be {\em transcribed/transformed} into a finite dimensional representation via discretization. For instance, differential equations and associated domains are often discretized using finite difference and quadrature schemes \cite{biegler2007overview}, while expectation operators and associated random domains are often discretized using Monte Carlo (MC) sampling and quadrature schemes \cite{chen2015scenario,kleywegt2002sample}. The finite-dimensional representation of the problem can be handled using standard optimization solvers (e.g., \texttt{Ipopt} and \texttt{Gurobi}). Sophisticated transformation techniques are used in different scientific disciplines; for example, projection onto orthogonal basis functions is a technique that is commonly used in PDE and stochastic optimization \cite{devolder2010solving, koivu2010galerkin}. 
\\

Although common mathematical elements of InfiniteOpt problems are found across disciplines, there are limited tools available to model and solve these problems in a unified manner.  Powerful, domain-specific software implementations are currently available for tackling dynamic and PDE optimization problems; examples include \texttt{Gekko}, \texttt{ACADO}, and \texttt{gPROMS} \cite{beal2018gekko, houska2011acado, asteasuain2001dynamic}. A key limitation of these tools is that the modeling abstraction used is specialized to specific problem classes and are not that easy to extend. On the other hand, there are powerful algebraic modeling languages such as \texttt{JuMP}, \texttt{CasADi}, \texttt{PETSc}, and \texttt{Pyomo} that offer high modeling flexibility to tackle different problem classes; however, these tools require the user to transform InfiniteOpt problems manually (which is tedious and prone to error). These limitations have recently motivated the development of modeling frameworks such as \texttt{Pyomo.dae} \cite{nicholson2018pyomo}; this framework unifies the notion of variables, objectives, and constraints defined over continuous, infinite-dimensional domains (sets).  This abstraction facilitates the modeling and transformation (via automatic discretization techniques) of optimization problems with embedded DAEs and PDEs. A limitation of this abstraction is that the notion of continuous sets is limited to space and time and this hinders generalizability (e.g., domains defined by random parameters need to be treated separately). Moreover, this framework provides limited capabilities to define measures (e.g., multi-dimensional integrals and risk functionals).  
\\

A unifying abstraction for InfiniteOpt  problems can facilitate the development of software tools and the development of new analysis and solution techniques. For instance, it has been recently shown that a graph abstraction unifies a wide range of optimization problems such as discrete-time dynamic optimization (graph is a line), network optimization (graph is the network itself), and multi-stage stochastic optimization (graph is a tree) \cite{jalving2019graph}. This unifying abstraction has helped identify new theoretical properties and decomposition algorithms \cite{jalving2020graph,shin2021exponential}; this has been enabled, in part,  via transferring techniques and concepts across disciplines.  The limited availability of unifying modeling tools ultimately limits our ability to experiment with techniques that appear in different disciplines and limits our ability to identify new modeling abstractions to tackle emerging applications. 
\vspace{0.2in}

\begin{figure}[!htb]
    \centering
    \includegraphics[width=0.9\textwidth]{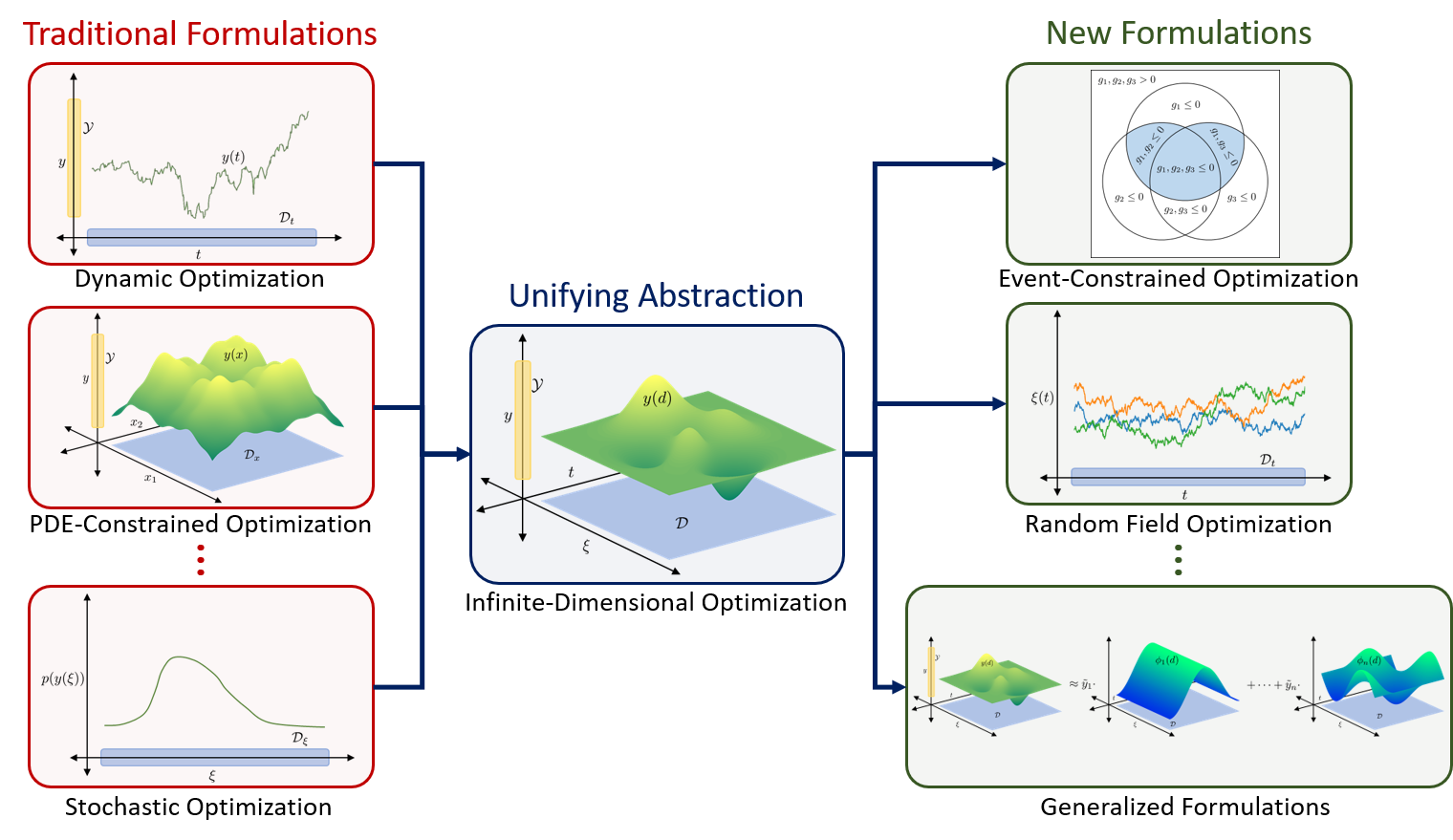}
    \caption{Summary of the proposed InfiniteOpt abstraction; the abstraction seeks to unify existing problem classes and use this to develop new classes.}
    \label{fig:abstract}
\end{figure}

In this work, we propose a unifying abstraction that facilitates the analysis, modeling, and solution of InfiniteOpt problems (see Figure \ref{fig:abstract}). Central to our abstraction is the notion of infinite parameters, which are parameters defined over infinite-dimensional domains (e.g., time, space, and random parameters). The proposed abstraction allows us to construct these domains in a systematic manner by using cartesian product operations and to define variables, objectives, and constraints over such domains and subdomains (restricted domains). The ability to handle restricted subdomains allows us to capture infinite-dimensional and standard (finite-dimensional) variables in a unified setting. Another key notion of the proposed abstraction are measure operators; these operators allow us to summarize infinite-dimensional objects over specific domains or subdomains and with this formulate problems with different types of objectives and constraints.  The proposed abstraction also incorporates differential operators, which are used to model how variables evolve other their corresponding domains. These modeling elements provide a bridge between different disciplines and enables cross-fertilization. For instance, we show that the proposed abstraction allows us to leverage the use of risk measures (used in stochastic optimization) to shape time-dependent trajectories (arising in dynamic optimization). The proposed abstraction also facilitates the development of new abstractions such as event-constrained optimization problems and optimization problems with space-time random fields. The proposed abstraction forms the basis of a {\tt Julia}-based modeling package that we call {\tt InfiniteOpt.jl}.  In this context, we show that the abstraction facilitates software development and enables a compact and intuitive modeling syntax and the implementation of transformation techniques (e.g., quadrature and sampling). 
\\

The paper is structured as follows. Section \ref{sec:abstraction} details the proposed unifying abstraction and highlights its implementation in \texttt{InfiniteOpt.jl}. Section \ref{sec:solution} reviews problem transformations into finite-dimensional representations through the lens of the abstraction. Section \ref{sec:consequences} discusses benefits provided by the abstraction. Section \ref{sec:cases} presents diverse case studies and Section \ref{sec:conclusion} closes the paper.  

\section{InfiniteOpt Abstraction} \label{sec:abstraction}

In this section, we outline the proposed unifying abstraction for InfiniteOpt problems. Specifically, we discuss the different elements of the abstraction, which include infinite domains and parameters, decision variables, measure operators, and differential operators. 

\subsection{Infinite Domains and Parameters} \label{sec:infinite_domains}

An InfiniteOpt problem includes a collection of infinite domains $\mathcal{D}_\ell$ with index $\ell \in \mathcal{L}$. An individual infinite domain is defined as $\mathcal{D}_\ell \subseteq \mathbb{R}^{n_\ell}$. The term {\em infinite} refers to the fact that an infinite domain has infinite cardinality (i.e., $|\mathcal{D}_\ell | = \infty$) and is thus continuous.  We also note that each infinite domain $\mathcal{D}_\ell$ is a subdomain of an $n_\ell$-dimensional Euclidean space $\mathbb{R}^{n_\ell}$.
\\

An infinite domain encompasses different domains encountered in applications; for instance, this can represent a time domain of the form $\mathcal{D}_t=[t_0, t_f]\subset\mathbb{R}$ (with $n_t=1$), a 3D spatial domain $\mathcal{D}_x=[-1, 1]^3\in \mathbb{R}^3$ (with $n_x=3$), or the co-domain of a multivariate random variable $\mathcal{D}_\xi=(-\infty,\infty)^m\in \mathbb{R}^{m}$ (with $n_\xi=m$).
\\

Infinite parameters are parameters that live in the associated infinite domains; specifically, we define the general parameter $d \in \mathcal{D}_\ell$. In our previous examples, the parameters are time $t \in \mathcal{D}_t$, space $x\in \mathcal{D}_x$, and random parameters $\xi \in \mathcal{D}_\xi$.
\\

\begin{figure}[!htb]
    \centering
    \includegraphics[width=0.7\textwidth]{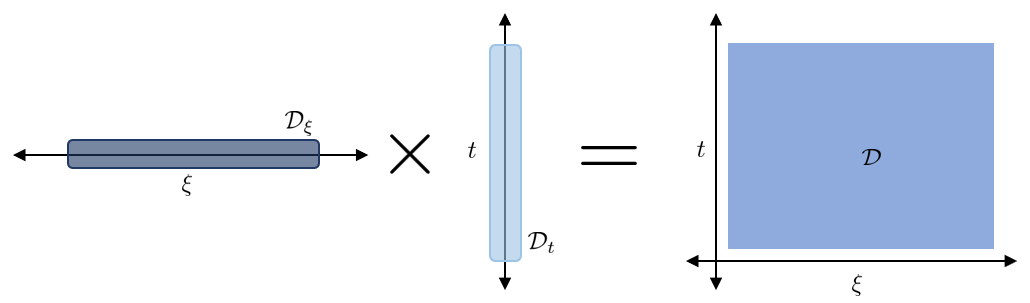}
    \caption{Cartesian product of random domain $\mathcal{D}_\xi$ and time domain $\mathcal{D}_t$ to produce $\mathcal{D}$.}
    \label{fig:infinite_domain}
\end{figure}

The domain of the InfiniteOpt problem is the cartesian product of the infinite domains:
\begin{equation}
    \mathcal{D} := \prod_{\ell \in \mathcal{L}} \mathcal{D}_\ell.
    \label{eq:infinite_domain}
\end{equation}
The construction of the problem domain is exemplified in Figure \ref{fig:infinite_domain}; here, we see that the domain obtained from the cartesian product of a 1D random domain $\mathcal{D}_\xi$ and a 1D temporal domain $\mathcal{D}_t$. 

\subsection{Decision Variables} \label{sec:variables}

A key feature of an InfiniteOpt problem is that it contains decision variables that are {\em functions} of infinite-dimensional parameters; as such, decision variables are also infinite-dimensional. In the proposed abstraction, we also consider InfiniteOpt problems that contain finite-dimensional variables (that do not depend on any parameters) and finite-dimensional variables that originate from reductions of infinite variables (e.g., integration over a domain or evaluation of an infinite variable at a point in the domain). To account for these situations, we define different types of variables: infinite, semi-infinite, point, and finite. 
\vspace{0.1in}

\begin{figure}[!htb]
    \centering
    \includegraphics[width=0.8\textwidth]{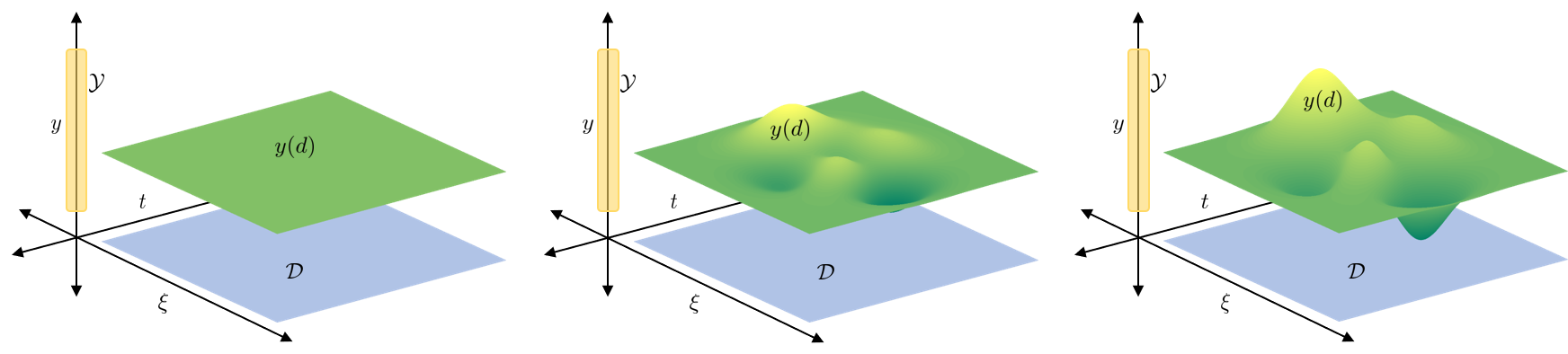}
    \caption{Depiction of realizations of an infinite variable $y(d)$ with $\mathcal{Y} \subseteq \mathbb{R}$ and $\mathcal{D} \subseteq \mathbb{R}^2$. The horizontal axes define the domain $\mathcal{D}$ and the vertical axis denotes the domain of feasible decisions $\mathcal{Y}$.}
    \label{fig:infinite_variable}
\end{figure}

\begin{figure}[!htb]
    \centering
    \begin{subfigure}[b]{0.3\textwidth}
        \centering
        \includegraphics[width=\textwidth]{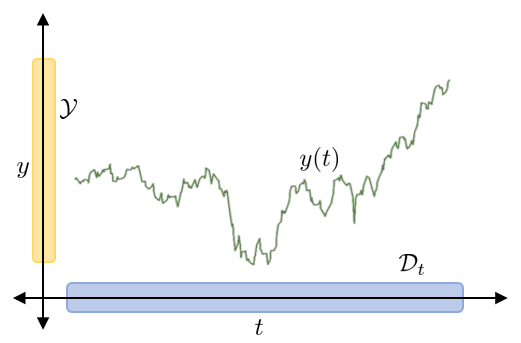}
        \caption{Dynamic Optimization}
        \label{fig:dynamic_variable}
    \end{subfigure}
    \quad
    \begin{subfigure}[b]{0.3\textwidth}
        \centering
        \includegraphics[width=\textwidth]{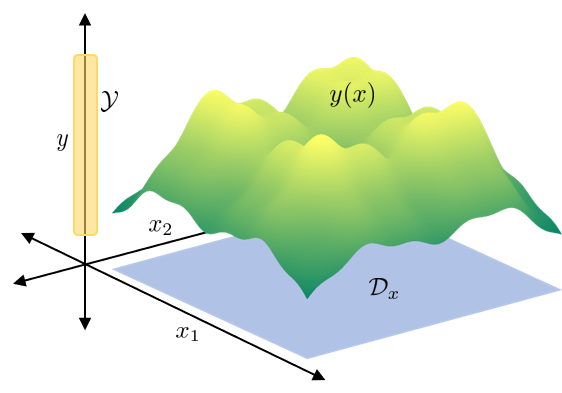}
        \caption{PDE Optimization}
        \label{fig:pde_variable}
    \end{subfigure}
    \quad
    \begin{subfigure}[b]{0.3\textwidth}
        \centering
        \includegraphics[width=\textwidth]{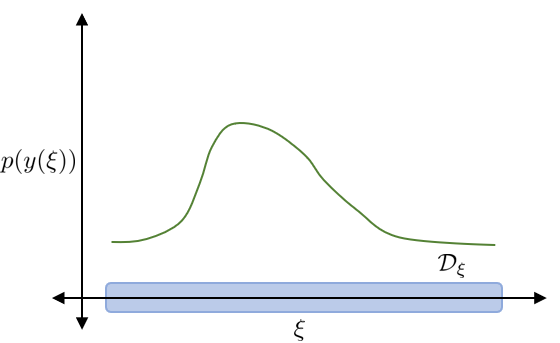}
        \caption{Stochastic Optimization}
        \label{fig:stochastic_variable}
    \end{subfigure}
    \caption{Example of infinite variables arising in traditional formulations. The recourse variable $y(\xi)$ is visualized in terms of its probability density function (as is customary in stochastic optimization).} 
    \label{fig:variable_examples}
\end{figure}

Infinite variables $y : \mathcal{D} \mapsto \mathcal{Y} \subseteq \mathbb{R}^{n_y}$ are functions that map an infinite domain $\mathcal{D}$ to the domain $\mathcal{Y}$. These variables are expressed as: 
\begin{equation}
    y(d) \in \mathcal{Y},\; d \in \mathcal{D}.
    \label{eq:infinite_variable}
\end{equation}
Figure \ref{fig:infinite_variable} shows that infinite variables are functions that can be interpreted as manifolds (also known as surfaces and fields). The goal of the InfiniteOpt problem is to shape these manifolds to achieve pre-determined goals (e.g., minimize their mean value or their peaks). Example classes of infinite variables include uncertainty-dependent decision policies arising in stochastic optimization (i.e., recourse variables), time-dependent control/state policies arising in dynamic optimization, or space-time fields arising in PDE optimization. For instance, in a stochastic PDE optimization problem, one might have an infinite variable of the form $y(t, x, \xi)$, which is simultaneously parameterized over time $t$, space $x$, and uncertainty $\xi$. Figure \ref{fig:variable_examples} illustrates some infinite variables commonly encountered in different disciplines.

\begin{figure}[!htb]
    \centering
    \begin{subfigure}[b]{0.3\textwidth}
        \centering
        \includegraphics[width=\textwidth]{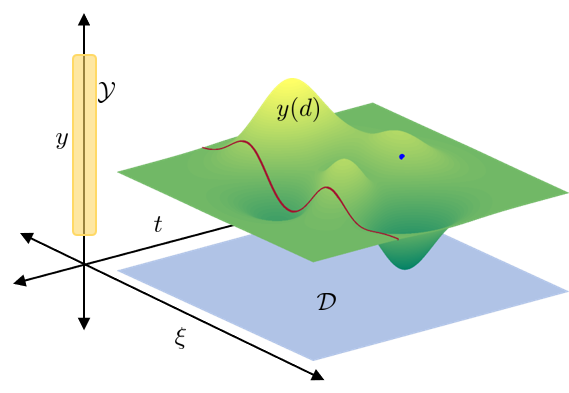}
        \caption{Infinite}
        \label{fig:infinite_variable_annotated}
    \end{subfigure}
    \quad \quad \quad
    \begin{subfigure}[b]{0.3\textwidth}
        \centering
        \includegraphics[width=\textwidth]{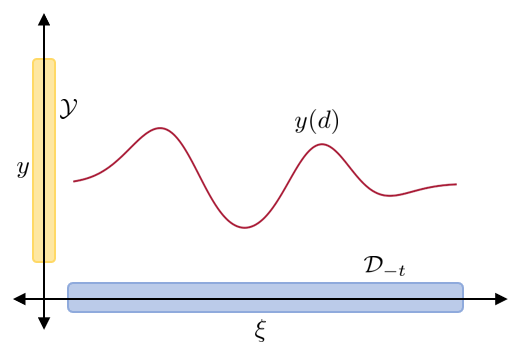}
        \caption{Semi-Infinite}
        \label{fig:semi-infinite_variable}
    \end{subfigure}
    \quad
    \begin{subfigure}[b]{0.182\textwidth}
        \centering
        \includegraphics[width=\textwidth]{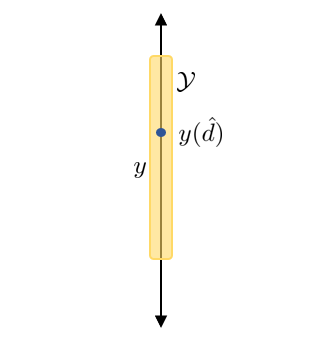}
        \caption{Point}
        \label{fig:point_variable}
    \end{subfigure}
    \caption{Illustration of how a semi-infinite variable $y(d),\; d\in \mathcal{D}_{-\ell}$ and a point variable $y(\hat{d})$ are obtained from an infinite variable $y(d),\, d\in \mathcal{D}$ via restriction/projection. Semi-infinite and point variables are realizations of an infinite variable and live in the domain $\mathcal{Y}$.}
    \label{fig:variable_summary}
\end{figure}

Semi-infinite variables $y : \mathcal{D}_{-\ell} \mapsto \mathcal{Y} \subseteq \mathbb{R}^{n_y}$ correspond to infinite variables in which the subdomain $\mathcal{D}_\ell$ has been restricted/projected to a single point $\hat{d}_\ell$; the restricted domain is denoted as $\mathcal{D}_{-\ell}$. These variables are also functions that map from the infinite domain $\mathcal{D}_{-\ell}$ to the domain $\mathcal{Y}$:
\begin{equation}
    y(d)  \in \mathcal{Y},\; d \in \mathcal{D}_{-\ell}.
    \label{eq:semi_infinite_variable}
\end{equation}
We refer to $\hat{d}_\ell \in \mathcal{D}_\ell$ as a support point of the domain. A depiction of how semi-infinite variables are derived from infinite variables via projection is provided in Figure \ref{fig:variable_summary}. Example classes mirror those of infinite variables with multiple parameter dependencies; for example, in a stochastic PDE problem, we might want to evaluate the variable $y(t,x,\xi)$ at the support point $t=0$ (initial time); this gives the semi-infinite variable $y(0, x, \xi)$ and  domain $\mathcal{D}_{-t}=\mathcal{D}_x\times \mathcal{D}_\xi$. 

\begin{figure}[!htb]
    \centering
    \includegraphics[width=0.3\textwidth]{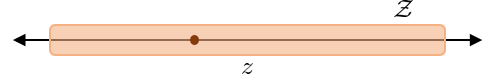}
    \caption{A depiction of a finite variable $z \in \mathcal{Z} \subseteq \mathbb{R}$.}
    \label{fig:finite_variable}
\end{figure}

Point variables denote infinite variables in which the entire $\mathcal{D}$ is restricted to a single point $\hat{d}$.  These are finite-dimensional variables that can be expressed as: 
\begin{equation}
    y(\hat{d})\in  \mathcal{Y}.
    \label{eq:point_variable}
\end{equation}
Figure \ref{fig:variable_summary} illustrates how point variables relate to other variable classes. Examples of point variables include random variables evaluated at a specific scenario/sample/realization of the uncertain parameter or a space-time variable evaluated at a specific point in the domain (e.g., boundary conditions). Point variables are also used in dynamic optimization to specify so-called point constraints (these ensure that a time trajectory satisfies a set of constraints at specific time points). Point variables can be thought of as variables that are held constant in an infinite domain. 
\\

The proposed abstraction also considers finite variables:
\begin{equation}
    z \in \mathcal{Z} \subseteq \mathbb{R}^{n_z}
    \label{eq:finite_variable}
\end{equation}
where $\mathcal{Z}$ denotes the feasible set. These are variables that are not parameterized over an infinite domain. Examples of finite variables arising in applications are first-stage decision variables and design variables arising in stochastic optimization or model parameters estimated in a dynamic optimization problem.  Figure \ref{fig:finite_variable} shows that this variable is analogous to the point variable depicted in Figure \ref{fig:point_variable}. In fact, this highlights that a point variable is a finite variable; the difference is that a point variable is derived from a restriction of an infinite variable.  However, we note that, technically speaking, a finite variable can be seen as a special case of an infinite variable in which the domain is a point. As such, infinite variables provide a unifying framework to capture different types of variables. 

\subsection{Differential Operators} \label{sec:derivatives}

Differential operators are a typical modeling element of InfiniteOpt problems. These operators capture how a given decision variable (an infinite-dimensional function) changes over its associated domain; for example, in a dynamic optimization problem, we want to know how quickly does a state/control variable change in time. 

\begin{figure}[!htb]
    \centering
    \includegraphics[width=0.7\textwidth]{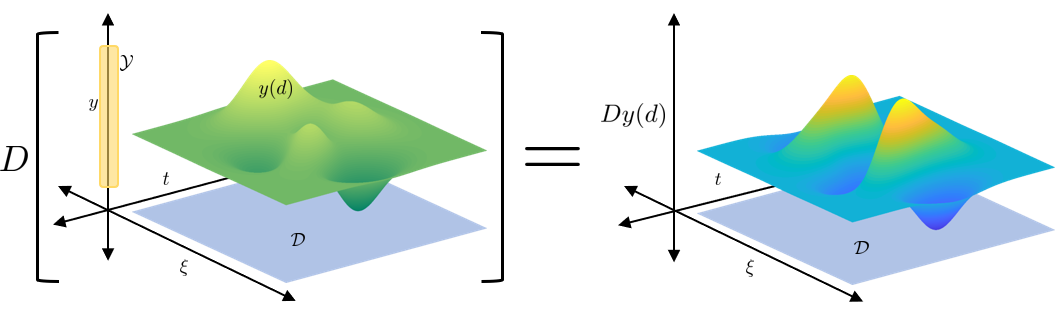}
    \caption{Depiction of a differential operator $D$ acting on the infinite variable $y(d)$.}
    \label{fig:derivative}
\end{figure}

The proposed abstraction defines a differential operator of the general form: 
\begin{equation}
    D: \mathcal{D} \mapsto \mathbb{R}. 
    \label{eq:derivative}
\end{equation}
Differential operators are applied to infinite variables $y(d),\; d\in \mathcal{D}$. These operators capture partial derivatives over individual parameters and more sophisticated operators that simultaneously operate on multiple parameters such as the Laplacian operator (typically encountered in PDE optimization). Figure \ref{fig:derivative} illustrates a differential operator acting on an infinite variable. 
\\

Differential operators map the infinite domain $\mathcal{D}$ to the scalar domain $\mathbb{R}$. The output of a differential operator is an infinite variable $D y(d),\; d\in \mathcal{D}$ that inherits the domain of the argument $y(d)$. For example, consider the infinite variable $y(t),\, t\in \mathcal{D}_t$; the partial derivative operator is an infinite variable $y'(t):=\partial y(t)/\partial t,\; t\in \mathcal{D}_t$. Some other specific examples include the partial derivative $\frac{\partial y(t, x, \xi)}{\partial t},\; (t,x,\xi)\in \mathcal{D}$ and the Laplacian $\Delta y(x),\; x \in \mathcal{D}_x$. Note that derivatives of the form $\frac{\partial}{\partial \xi}$ are not typically used in stochastic optimization problems; however, the proposed abstraction allows for this operator to be defined. This modeling feature can be used, for instance, to control how a random variable changes in the uncertainty space (this can be used to manipulate the shape of its probability density function). 

\subsection{Measure Operators} \label{sec:measures}
\begin{figure}[!htb]
    \centering
    \includegraphics[width=0.75\textwidth]{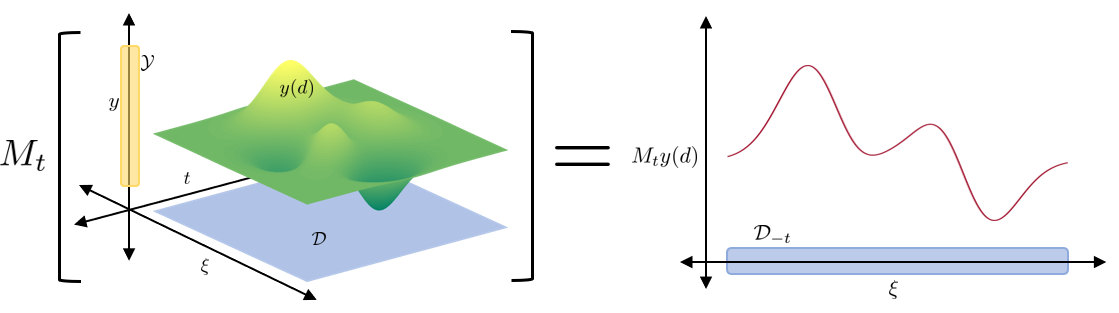}
    \caption{Measure operator $M_t$ that acts on domain $\mathcal{D}_t$ of the infinite variable $y(d),\; d\in \mathcal{D}$. This operation returns the semi-infinite variable $m(\xi):=M_t\,y,\; \xi \in \mathcal{D}_\xi$.}
    \label{fig:measure}
\end{figure}

Measure operators are key modeling constructs that are used to {\em summarize} functions by collapsing them to a single quantity. For example, in a dynamic optimization problem, one typically minimizes the time-integral of the cost (a scalar quantity). The proposed abstraction defines a measure operator of the general form: 
\begin{equation}
    M_\ell: \mathcal{D} \mapsto \mathbb{R}
    \label{eq:measure}
\end{equation}
Here, the index $\ell$ indicates that the operator is applied on the subdomain $\mathcal{D}_\ell$ and thus has the effect of restricting the domain. As such, the output of a measure operator is a semi-infinite variable that lives in the restricted domain $\mathcal{D}_{-\ell}$. Figure \ref{fig:measure} illustrates such a measure operator.
\vspace{0.2in}

Measure operators are a key feature of InfiniteOpt problems; specifically, objective functions and constraints are often expressed in terms of measure operators. For instance, consider a field $y(t,x,\xi)$ arising in a stochastic PDE problem; one can define measure operator that computes the time-integral $m(x,\xi):=\int_{t \in \mathcal{D}_t} y(t, x, \xi) dt,\; (x, \xi) \in \mathcal{D}_x \times \mathcal{D}_\xi$ and note that the output of this operation is a semi-infinite variable $m(x,\xi)$ that lives in $\mathcal{D}_{-t}=\mathcal{D}_{x}\times \mathcal{D}_\xi$. One can also define an operator that computes the expectation $m(t,x):=\int_{\xi \in \mathcal{D}_\xi} y(t, x, \xi) p(\xi)d\xi,\; (t,x) \in \mathcal{D}_t \times \mathcal{D}_x$ (where $p(\cdot)$ is the probability density function of $\xi$);  this operation gives a semi-infinite variable $m(t,x)$ that lives in $\mathcal{D}_{-\xi}=\mathcal{D}_t\times \mathcal{D}_x$.
\\

The expectation is a measure operator that is of particular interest in stochastic optimization because this can be used to compute different types of risk measures and probabilities; for instance, in the previous example, one might want to compute a probability of the form:
\begin{equation}
    \mathbb{P}_\xi(y(t, x, \xi) \in \mathcal{Y}),\;\; (t, x) \in \mathcal{D}_t \times \mathcal{D}_x.
\end{equation}
This is the probability that $y(t, x, \xi)$ is in the domain $\mathcal{Y}$ and can be computed by using an expectation operator:
\begin{equation}
    \begin{aligned}
        \mathbb{P}_\xi(y(t, x, \xi) \in \mathcal{Y})&=\mathbb{E}_\xi\big[\mathbbm{1}[y(t, x, \xi)\in \mathcal{Y}] \big]\\
                                                    &=\int_{\xi \in \mathcal{D}_\xi} \mathbbm{1}[y(t, x, \xi)\in \mathcal{Y}] p(\xi)d\xi.
    \end{aligned}
\end{equation}
where $\mathbbm{1}[\cdot]$ is the indicator function  and the argument of this function is the {\em event} of interest.  We recall that the indicator function returns a value of 0 if the event is not satisfied or a value of 1 if the event is satisfied. An important observation is that the indicator function can be used to define a wide range of measures and over different types of domains; for instance, the measure $\int_{t \in \mathcal{D}_t} \mathbbm{1}[y(t)> \overline{y}] dt$ denotes the amount of time that the function $y(t),\; t\in \mathcal{D}_t$ crosses the threshold $\bar{y}$.
\\

Measure operators can also be used to summarize infinite variables over multiple subdomains; for example, one can consider the following measures:
\begin{subequations}
    \begin{align}
        M_{t,x}\,y&=\int_{t \in \mathcal{D}_t}\int_{x \in \mathcal{D}_x}y(t, x, \xi) dx dt, \;\; \xi \in \mathcal{D}_\xi\\
        M_{t,x,\xi}\, y&=\mathbb{E}_{\xi}\left[\int_{t \in \mathcal{D}_t}\int_{x \in \mathcal{D}_x}y(t, x, \xi) dx dt\right]\\
        M_{t,x,\xi}\, y&=\mathbb{E}_{\xi}\left[\int_{t \in \mathcal{D}_t}\int_{x \in \mathcal{D}_x} \mathbbm{1}[y(t, x, \xi)\in \mathcal{Y}]  dx dt\right].
    \end{align}  
\end{subequations}
One can thus see that a wide range of measures can be envisioned. 

\subsection{Objectives} \label{sec:objective}

In InfiniteOpt problems, objective functions are functions of infinite variables; as such, objectives are infinite variables  (infinite-dimensional functions). Minimizing or maximizing an infinite-dimensional function does not yield a well-posed optimization problem. This situation is similar in spirit to that appearing in multi-objective optimization problem, in which we seek to simultaneously minimize/maximize a finite collection of objectives (in an InfiniteOpt problem, the collection is infinite).
\\

To deal with ill-posedness, one often resorts to scalarization techniques; the idea is to reduce/summarize the infinite-dimensional function into a single scalar quantity. The goal of this scalarization procedure is to manipulate the shape of the infinite-dimensional objective (e.g., minimize its mean value or its extreme value). Scalarization is performed by using measure operators; for instance, in the context of multi-objective optimization, one scalarizes the objectives by computing a weighted summation of the objectives. In an InfiniteOpt setting, this weighting is done by computing a weighted integral of the objective. For instance, in dynamic optimization, we often have a time-dependent objective function $f(t):=f(y(t),t),\; t\in \mathcal{D}_t$; here, we can notice that the objective depends on an infinite variable and is thus also an infinite variable. We can scalarize this variable by using the measure $M_t f:=\int_{t \in \mathcal{D}_t} f(t)w(t)dt$ with a weighting function satisfying $w: \mathcal{D}_t \to  [0,1]$ and $\int_{t\in \mathcal{D}_t} w(t)dt=1$ (note that this measure is a time-average of the objective trajectory). 
\\

In space-time PDE optimization, the objective is defined over an infinite domain $\mathcal{D}=\mathcal{D}_{t}\times \mathcal{D}_x$ that depends on decision variables $y(t,x) \in \mathcal{Y}$; as such, the objective is given by the field $f(t,x):=f(y(t,x),t,x),\; t\in \mathcal{D}_t, x\in \mathcal{D}_x$. One can scalarize this field by using a measure:
\begin{equation}
    M_{t, x} f = \int_{(t,x)\in\mathcal{D}_{t,x}} f(t,x) w(t,x) dt dx,
    \label{eq:optimal_control_obj}
\end{equation}
with weighting function satisfying $w: \mathcal{D}_{t,x}\to [0,1]$ and $\int_{(t, x) \in \mathcal{D}_{t,x}}w(t,x)dtdx=1$. One can think of this measure as a space-time average of the objective.
\\

\begin{figure}[!htb]
    \centering
    \includegraphics[width=0.55\textwidth]{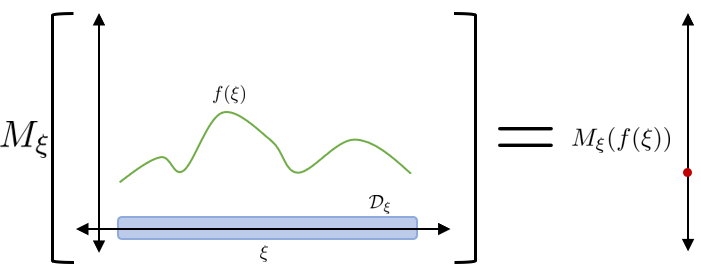}
    \caption{Depiction of measure operator $M_\xi$ acting on infinite variable $f(\xi) = f(z,y(\xi))$.}
    \label{fig:expectation}
\end{figure}

In stochastic optimization, we have infinite-dimensional objectives of the form $f(\xi) := f(z,y(\xi)),$ $\xi \in \mathcal{D}_\xi$, where $y(\xi)$ is a recourse variable (an infinite variable). Scalarization can be achieved by using the expectation operator:
\begin{equation}
    M_\xi f = \mathbb{E}_\xi[f(\xi)] = \int_{\xi \in \mathcal{D}_\xi} f(\xi) p(\xi) d\xi
    \label{eq:expectation}
\end{equation}
where $p(\xi)$ is the probability density function satisfying $p(\xi)\geq 0$ and $\int_{\xi \in \mathcal{D}_\xi}p(\xi)=1$.  This measure is illustrated in Figure \ref{fig:expectation}.

Average measures as those described previously are intuitive and widely used in practice; however, in Section \ref{sec:generalized_measures} we will see that one can envision using a huge number of measures to perform scalarization. The huge number of choices arises from the fact that one can manipulate the shape of an infinite-dimensional function in many different ways (by focusing in different features of the function); for instance, one can minimize the peak of the function or minimize its variability.  In the field of stochastic optimization,  for instance, one aims to manipulate the shape of the infinite-dimensional objective by selecting different risk measures (summarizing statistics) such as the variance, median, quantile, worst/best -case value, or probabilities. We will see that one can borrow risk measures used in stochastic optimization to summarize infinite variables in other domains (e.g., space-time); this leads to interesting approaches to shape complex manifolds/fields arising in complex InfiniteOpt problems. 

\subsection{Constraints} \label{sec:constrs}

\begin{figure}[!htb]
    \centering
    \includegraphics[width=0.4\textwidth]{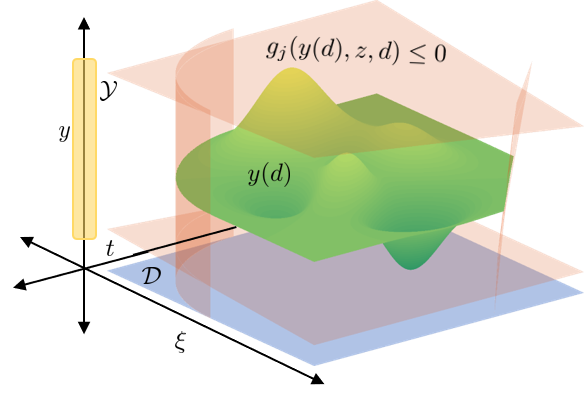}
    \caption{Depiction of infinite-dimensional constraints $g_j(y(d), z, d) \leq 0$ defined over an infinite domain $\mathcal{D}$.}
    \label{fig:constraints}
\end{figure}

As in the case of objectives, constraints in InfiniteOpt  problems depend on infinite variables and are thus infinite variables themselves.  One thus need to use specialized techniques to handle constraints and with this ensure that the problem is well-posed. A key observation that arises in this context is that constraints are treated differently than objectives; specifically, one typically seeks to impose bounds on constraint values and one can handle collections of constraints simultaneously. For instance, in semi-infinite optimization problems, one enforces constraints of the form: 
\begin{equation}
    g_j(y(d), d) \leq 0,\; j \in \mathcal{J}, d \in \mathcal{D}.
    \label{eq:semi_inf_constrs}
\end{equation}
In vector form, this collection of constraints can be expressed as:
\begin{equation}
    g(y(d), d) \leq 0,\; d \in \mathcal{D}.
    \label{eq:semi_inf_constrs_vect}
\end{equation}
where $g(\cdot)$ is a vector function that contains the constraint collection $g_j(\cdot)\, j\in \mathcal{J}$. We can see that the constraint functions $g(\cdot)$ are required to a take value below zero {\em for all} values of the parameter $d \in \mathcal{D}$.  Moreover, we can see that the constraints $j\in \mathcal{J}$ are all enforced at once. Figure \ref{fig:constraints} illustrates this constraint system. This particular approach to enforcing constraints is also widely used in dynamic optimization and stochastic optimization. For instance, in the context of dynamic optimization, one may seek to keep time trajectories for controls/states below a certain threshold value for all times in a time horizon. In the context of stochastic optimization, one may seek to satisfy the demand of a product for all possible realizations of uncertainty (in this context the constraints are said to be enforced {\em almost surely} or with probability of one).  

These types of constraints are defined in our abstraction using the general form:
\begin{equation}
    g(Dy, y(d), z, d) \leq 0, \; d \in \mathcal{D}.
    \label{eq:domain_constrs}
\end{equation}
These encapsulate the above use cases and are exemplified by the following PDE optimization constraints that include differential operators, path constraints, and point constraints (e.g., boundary conditions):
\begin{equation}
    \begin{gathered}
        g\left(Dy(t, x), y(t, x), t, x\right) = 0, \; (t, x) \in \mathcal{D}_t \times \mathcal{D}_x \\
        g\left(y(t, x), t, x\right) \leq 0, \; (t, x) \in \mathcal{D}_t \times \mathcal{D}_x \\
        g\left(y(\hat{t}, \hat{x}), \hat{t}, \hat{x}\right) \leq 0. 
    \end{gathered}
    \label{eq:space-time_constrs}
\end{equation}

Constraints that follow the form of \eqref{eq:domain_constrs} can be quite restrictive for certain applications, since they need to hold for every value parameter $d \in \mathcal{D}$. One can relax this requirement by instead enforcing the constraint on a selected set of points in the domain $\mathcal{D}$ or by enforcing constraints on a measure of the constraint functions. For instance, consider the set of constraint functions $h_k(Dy(d), y(d), z, d), \; k \in \mathcal{K}$; we can aim to enforce  constraints on expected values of such functions as:
\begin{equation}
    \int_{\xi \in \mathcal{D}_\xi} \; h_k(y(\xi), \xi) p(\xi)d\xi \geq 0, \; k \in \mathcal{K}.
    \label{eq:expect_constrs}
\end{equation}
Given that there are a wide range of measures that can help shape infinite-dimensional functions, one can also envision different approaches to enforce constraints. For instance, in stochastic optimization, one typically uses scalar chance (probabilistic) constraints of the form: 
\begin{equation}
    \mathbb{P}_\xi\left (h_k(y(\xi), \xi) \leq 0\right) \geq \alpha, \; k \in \mathcal{K}. 
    \label{eq:single_chance_constrs}
\end{equation}
This set of constraints require that each constraint function $h_k(\cdot)$ is kept below zero to a certain probability level $\alpha$. In stochastic optimization, one also often enforces joint chance constraints:
\begin{equation}
    \mathbb{P}_\xi\left (h_k(y(\xi), \xi) \leq 0, \;  k \in \mathcal{K}\right) \geq \alpha.
    \label{eq:joint_chance_constr1}
\end{equation}
The joint chance constraint can also be expressed in vector form as:
\begin{equation}
    \mathbb{P}_\xi\left (h(y(\xi), \xi) \leq 0\right) \geq \alpha.
    \label{eq:joint_chance_constr2}
\end{equation}
Joint chance constraints require that the constraint functions $h(\cdot)$ are kept (jointly) below zero with a certain probability level $\alpha$. We will see that joint chance constraints allow us to enforce constraints on probability of {\em events} and we will see that this provides a flexible modeling construct to capture complex decision-making logic. For instance, we might want to ensure that the temperature of a system is higher than a certain value {\em and} that the concentration of the system is lower than a certain value with a given probability. Joint chance constraints can also be interpreted as a generalization of other constraint types; for instance, if we set $\alpha=1$, the constraint \eqref{eq:joint_chance_constr1} is equivalent to \eqref{eq:semi_inf_constrs}.
\\

The above measure constraints can be expressed in the following general form:
\begin{equation}
    M h(Dy, y(d), z, d) \geq 0.
    \label{eq:meas_constrs}
\end{equation}
For instance, the chance constraint \eqref{eq:joint_chance_constr2} can be expressed as:
\begin{equation}
    M_\xi h \geq 0
    \label{eq:chance_constr}
\end{equation}
with
\begin{equation}
    M_\xi h = \mathbb{E}_\xi\left[\mathbbm{1}[h(y(\xi), \xi) \leq 0]\right] - \alpha. 
\end{equation}

\subsection{InfiniteOpt Formulation} \label{sec:infinite_problem}

We summarize the previous elements to express the InfiniteOpt problem in the following abstract form: 
\begin{equation}
    \begin{aligned}
        &&\min_{y(\cdot)\in\mathcal{Y}, z\in \mathcal{Z}} &&& M f\left(Dy, y(d), z, d\right) \\
        && \text{s.t.} &&& g\left(Dy, y(d), z, d\right) \leq 0, \; d \in \mathcal{D} \\
        &&&&& M h\left(Dy, y(d), z, d\right) \geq 0. 
    \end{aligned}
    \label{eq:infinite_form}
\end{equation}
This abstract form seeks to highlight the different elements of the proposed abstraction (e.g., infinite domains and variables, finite variables, measure operators, differential operators).

\subsection{Implementation in \texttt{InfiniteOpt.jl}} \label{sec:infiniteopt_modeling}

We now proceed to describe how the proposed abstraction can facilitate the development of modeling tools. Specifically, the proposed abstraction is used as the backbone of a modeling package that we call \texttt{InfiniteOpt.jl} (\url{https://github.com/zavalab/InfiniteOpt.jl}). \texttt{InfiniteOpt.jl} is written in the {\tt Julia}  programming language \cite{bezanson2017julia} and builds upon the capabilities of {\tt JuMP.jl} \cite{dunning2017jump} to intuitively and compactly express InfiniteOpt problems.
\\

Some of the modeling features of \texttt{InfiniteOpt.jl} are illustrated by using the example problem:
\begin{subequations}
    \begin{align}
        &&\min_{y_a(t), y_b(t,\xi), y_c(\xi), z} &&& \int_{t \in \mathcal{D}_t} y_a(t)^2 + 2 \mathbb{E}_{\xi}[y_b(t, \xi)] dt \label{eq:ex_obj}\\
        && \text{s.t.} &&& \frac{\partial y_b(t, \xi)}{\partial t} = y_b(t,\xi)^2 + y_a(t) - z_1, && \forall t \in \mathcal{D}_t, \ \xi \in \mathcal{D}_\xi \label{eq:ex_constr1}\\
        &&&&& y_b(t, \xi) \leq y_c(\xi)U, && \forall t \in \mathcal{D}_t \\
        &&&&& \mathbb{E}_\xi[y_c(\xi)] \geq \alpha \\
        &&&&& y_a(0) + z_2 = \beta \label{eq:ex_constr4}\\
        &&&&& y_a(t),y_b(t, \xi) \in \mathbb{R}_+, y_c(\xi) \in \{0, 1\}, z \in \mathbb{Z}^2, && \forall t \in \mathcal{D}_t, \ \xi \in \mathcal{D}_\xi \label{eq:ex_var_props}
    \end{align}
    \label{eq:example_model}
\end{subequations}
Here, $y_a(t)$, $y_b(t, \xi)$, and $y_c(\xi)$ are infinite variables, $z$ are finite variables, $U,\alpha, \beta \in \mathbb{R}$ are constants, $\mathcal{D}_t=[t_0,t_f]$ is the time domain, and $\mathcal{D}_\xi$ is the co-domain of the random parameter $\mathcal{N}(\mu, \Sigma)$.
\\

The corresponding \texttt{InfiniteOpt.jl} syntax for expressing this problem is shown in Code Snippet \ref{code:infiniteopt_example}. An InfiniteOpt problem is stored in an  \texttt{InfiniteModel} object; line \ref{line:model_define} shows the initialization of the model object \texttt{model}. The model is automatically transcribed into a finite dimensional representation and solved using the \texttt{KNITRO} solver \cite{nocedal2006knitro}. More information on how the \texttt{InfiniteModel} is transcribed by \texttt{InfiniteOpt.jl} is provided in Section \ref{sec:infiniteopt_solution}. Lines \ref{line:define_t} and \ref{line:define_xi} use \texttt{@infinite\_parameter} to define the infinite parameters with their respective infinite domains and indicate that each domain should use 100 finite supports in the transcription. The random parameters $\xi$ can be associated with any probably density function supported by the Julia package \texttt{Distributions.jl} \cite{besanccon2019distributions}. Lines \ref{line:define_ya}-\ref{line:define_z} define the decision variables and their associated properties in accordance with Equation \eqref{eq:ex_var_props} following a symbolic \texttt{JuMP.jl}-like syntax by means of \texttt{@variable}. Line \ref{line:define_obj} defines the complex objective depicted in Equation \eqref{eq:ex_obj} via \texttt{@objective}. Lines \ref{line:define_c1}-\ref{line:define_c4} define constraints \eqref{eq:ex_constr1}-\eqref{eq:ex_constr4} using \texttt{@constraint}. Notice how the differential operator and measure operators (in this case an expectation and an integral) are easily incorporated using \texttt{Julia} syntax. Lines \ref{line:optimize}-\ref{line:get_obj} illustrate how the model \texttt{model} is solved using \texttt{optimize!} and then how the solution information can be extracted from the model.

\begin{minipage}[!htb]{0.9\linewidth}
\begin{scriptsize}
\lstset{language=Julia,breaklines = true}
\begin{lstlisting}[label = {code:infiniteopt_example},caption = Modeling problem \eqref{eq:example_model} using \texttt{InfiniteOpt.jl}.]
using InfiniteOpt, Distributions, KNITRO

# Initialize the model
model = InfiniteModel(KNITRO.Optimizer) |\label{line:model_define}|

# Add the infinite parameters corresponding to the infinite domains
@infinite_parameter(model, t ∈ [t0, tf], num_supports = 100) |\label{line:define_t}|
@infinite_parameter(model, ξ[1:10] ~ MvNormal(μ, Σ), num_supports = 100) |\label{line:define_xi}|

# Add the variables and their domain constraints
@variable(model, 0 ≤ ya, Infinite(t)) |\label{line:define_ya}|
@variable(model, 0 ≤ yb, Infinite(t, ξ))
@variable(model, yc, Infinite(ξ), Bin)
@variable(model, z[1:2], Int) |\label{line:define_z}|

# Define the objective 
@objective(model, Min, ∫(ya ^ 2 + 2 * 𝔼(yb, ξ), t)) |\label{line:define_obj}|

# Add the constraints 
@constraint(model, ∂(yb, t) == yb ^ 2 + ya - z[1]) |\label{line:define_c1}|
@constraint(model, yb ≤ yc * U)
@constraint(model, 𝔼(yc, ξ) ≥ α)
@constraint(model, ya(0) + z[2] == β) |\label{line:define_c4}|

# Solve and retrieve the results
optimize!(model) |\label{line:optimize}|
opt_objective = objective_value(model) |\label{line:get_obj}|
\end{lstlisting}
\end{scriptsize}
\end{minipage}

\section{InfiniteOpt Transformations} \label{sec:solution}

We now discuss how InfiniteOpt problems are solved through the lens of the proposed unifying abstraction. Solution approaches typically rely on transforming the InfiniteOpt problem \eqref{eq:infinite_form} into a finite-dimensional formulation that can be solved using conventional optimization solvers. There are a large number of possible methods to transform InfiniteOpt problems that are used in different domains such as dynamic, PDE, and stochastic optimization.  Our goal here is not to provide an extensive discussion and implementation of all these approaches; instead, we highlight common elements of different approaches and motivate how these can be facilitated by using a unifying abstraction. 

\subsection{Direct Transcription} \label{sec:transcription}

Our first goal is to obtain a finite representation of an infinite domain $\mathcal{D}_\ell$; direct transcription accomplishes this via a finite set of support points that we represent as $\hat{\mathcal{D}}_\ell = \{\hat{d}_{\ell, i} : \hat{d}_{\ell, i} \in \mathcal{D}_\ell,  i \in \mathcal{I}_\ell\}$ . The concept of the support set $\hat{\mathcal{D}}_\ell$ used here is general and a variety of methods can be employed to generate it. In stochastic optimization, for instance, a set of MC samples is typically drawn from a probability density function of the infinite parameters \cite{birge2011introduction}, while PDE problems commonly use quadrature schemes \cite{shin2020diffusing}. The proposed abstraction seeks to enable porting techniques across fields; for instance, one might envision generating support points for a space-time domain by sampling or one might envision generating support points for a random domain by using quadrature points (as done in sparse grids and Latin hypercube sampling).  

\begin{figure}[!htb]
    \centering
    \includegraphics[width=0.7\textwidth]{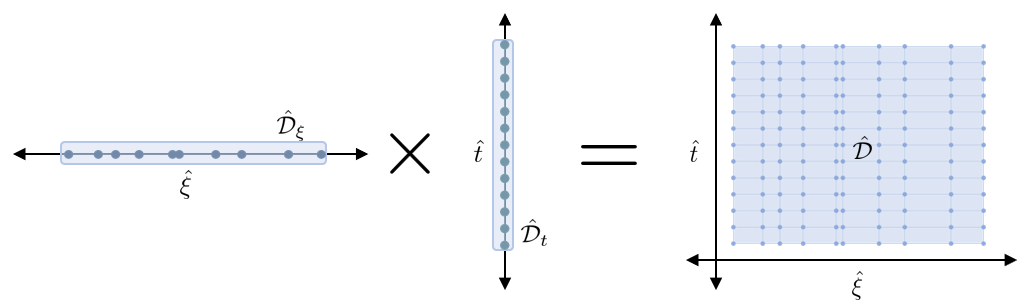}
    \caption{Finite support set $\hat{\mathcal{D}}$ that approximates infinite domain $\mathcal{D}$.}
    \label{fig:discrete_domain}
\end{figure}

The support set for the infiniteOpt problem $\hat{\mathcal{D}}$ is the cartesian product of the individual supports sets:
\begin{equation}
    \hat{\mathcal{D}} := \prod_{\ell \in \mathcal{L}} \hat{\mathcal{D}}_\ell.
    \label{eq:support_set}
\end{equation}
Figure \ref{fig:discrete_domain} illustrates how the support set (a finite domain) approximates the infinite domain $\mathcal{D}$. Note that this definition of $\hat{\mathcal{D}}$ assumes that the individual domains $\mathcal{D}_\ell$ are independent of one another. This assumption does not hold in some complex applications; for instance, in stochastic dynamic optimization problems, we might have random parameters that are functions of time (this is discussed further in Section \ref{sec:random_fields}). 

\begin{figure}[!htb]
    \centering
    \includegraphics[width=0.75\textwidth]{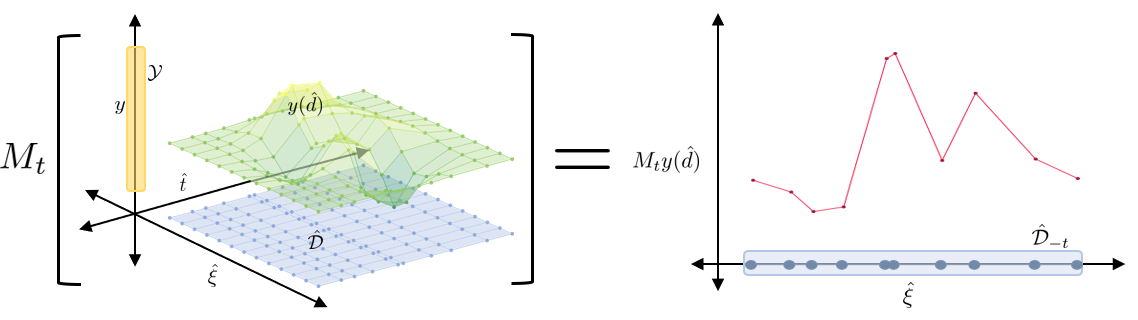}
    \caption{Depiction of a measure operator $M_t$ approximated via a numerical scheme (quadrature in this case).}
    \label{fig:discrete_measure}
\end{figure}

The infinite-dimensional formulation \eqref{eq:infinite_form} is projected onto the finite support set to yield a finite-dimensional approximation that can be modeled using conventional optimization solvers. We now proceed to discuss how this projection is achieved. Measures are approximated with an appropriate numerical scheme; this can take on a range of forms and may require the incorporation of additional supports and variables. For instance, a common set of measures (e.g., expectations and integrals over space-time domains) are of the form:
\begin{equation}
    M_\ell y = \int_{d' \in \mathcal{D}_\ell} y(d') w(d') dd'
    \label{eq:integral_measure}
\end{equation}
where $w(\cdot)$ is a weighting function. Such measures can be approximated using support points as:
\begin{equation}
    M_\ell y \approx \sum_{i \in \mathcal{I}_\ell} \beta_i y(\hat{d}_{\ell, i}) w(\hat{d}_{\ell, i}).
    \label{eq:integral_measure_finite}
\end{equation}
This general form is used in quadrature and sampling schemes; the only difference between these schemes arises in how the supports $\hat{d}_{\ell, i}$ and the coefficients $\beta_i$ are selected. Figure \ref{fig:discrete_measure} depicts a measure approximated via quadrature.

\begin{figure}[!htb]
    \centering
    \includegraphics[width=0.75\textwidth]{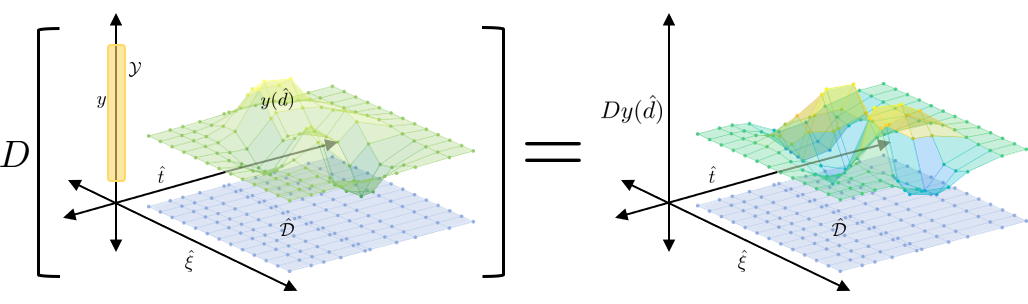}
    \caption{Depiction of a differential operator $D$ approximated via a numerical scheme (central finite differences in this case) relative to a realization of infinite variable $y(d)$.}
    \label{fig:discrete_derivative}
\end{figure}

Differential operators appearing in formulation \eqref{eq:infinite_form} also need to be approximated. Sometimes these operators can be reformulated in integral form; in such a case, one can use the measure approximations previously discussed. However, in some cases, this reformulation is not possible; 
for instance, a differential operator might be implicitly defined within expression functions (e.g., boundary conditions) and/or within measures. In our framework, we treat differential operators as infinite variables and handle them via {\em lifting}. To illustrate how this is done, suppose that we have an expression of the form:
\begin{equation}
g\left(\frac{\partial y(d)}{\partial d_\ell}, y(d), z\right)=0,\; d \in \mathcal{D}_\ell.
\end{equation}
Here, we introduce an auxiliary variable $y'(d)$ and reformulate the above constraint as:
\begin{subequations}
\begin{align}
&g(y'(d), y(d), z)=0,\; d \in \mathcal{D}_\ell\\
&\frac{\partial y(d)}{\partial d_\ell} = y'(d), \;d \in \mathcal{D}_\ell. 
\end{align}
\end{subequations}
The second expression can now be approximated using traditional schemes using support points; for instance, when $d$ denotes time (e.g., in a dynamic optimization problem), one typically uses a  backward finite difference:
\begin{equation}
    y(\hat{d}_{\ell, i}) = y(\hat{d}_{\ell, i-1}) + (\hat{d}_{\ell, i} - \hat{d}_{\ell, i-1})    y'(d_{\ell,i}). \label{eq:implicit_euler}
\end{equation}
Figure \ref{fig:discrete_derivative} illustrates how these techniques approximate differential operators. A lifting approach can be used to handle higher-order and multi-dimensional operators (e.g., Laplacian) via nested recursions.  These basic constructs can be used to enable the implementation of direct transcription schemes such as MC sampling, quadrature, finite difference, and orthogonal collocation over finite elements. 

Once the measures and derivatives are approximated, the direct transcription procedure is finalized by projecting the remaining constraints with infinite domain dependencies over the finite support set $\hat{\mathcal{D}}$. The transformation incurred by direct transcription is often linear since the typical measure and differential operator approximations are linear transformations of the respective modeling objects (e.g., MC sampling and finite difference). For instance, this means that if the InfiniteOpt problem of interest is an infinite quadratic program (QP), then its transcribed variant will typically be a finite QP. 

\subsection{Alternative Transformations} \label{sec:mwr}

Direct transcription is a common class of methods for transforming an InfiniteOpt problem into a finite-dimensional representation by using a finite set of support points.  A limitation of this approach is that it does not provide a solution in functional form (it only provides a solution defined at the support points).  Alternative transformation methods can be envisioned to deliver solutions in functional form. The method of weighted residuals (MWR) is a general class of methods that conducts the transformation by approximating the problem elements using basis expansions. Popular MWR techniques include polynomial chaos expansion (PCE) used in stochastic optimization \cite{xiu2010numerical} and orthogonal collocation used in dynamic optimization \cite{armaou2002dynamic,koivu2010galerkin}. For instance, Gnegel et. al. recently demonstrated how such basis expansion techniques can enhance the tractability of mixed-integer PDE problems relative to using traditional transcription methods \cite{gnegel2021solution}.

\begin{figure}[!htb]
    \centering
    \includegraphics[width=1\textwidth]{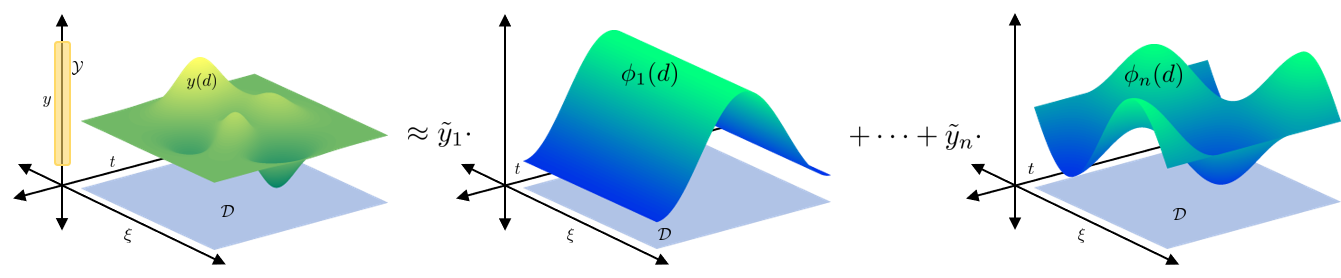}
    \caption{Depiction of how an infinite variable $y(d)$ can be approximated as a linear combination of basis functions $\phi_i(\cdot)$.}
    \label{fig:basis_functions}
\end{figure}

In MWR, a set of trial/basis functions $\Phi = \{\phi_i(d): i \in \mathcal{I}\}$ is defined over an infinite domain $\mathcal{D}$ and linear combinations of these functions are used to approximate the infinite variables:
\begin{equation}
    y(d) \approx \sum_{i \in \mathcal{I}} \tilde{y}_i \phi_i(d)
    \label{eq:function_basis}
\end{equation}
where $\tilde{y}_i \in \mathbb{R}$ are the basis function coefficients. An illustration of this approximation is given in Figure \ref{fig:basis_functions}; here, we require that the basis functions $\phi_i(d)$ and the infinite variables $y(d)$ both reside in a common space such that this approximation becomes exact when the set $\Phi$ is an orthogonal set of basis functions and $|\Phi| \rightarrow \infty$ \cite{graham2013modeling}. Since the basis functions are known, this representation allows us to represent the infinite variables $y(d)$ in terms of the coefficients $\tilde{y}_i$ (which are finite variables). As such, this approach effectively transforms infinite variables into finite variables. The goal is now to project the formulation \eqref{eq:infinite_form} onto a set of basis functions so as to obtain a finite formulation that solely depends on the finite variables $\tilde{y}_i$ and $z$. This is done by expressing differential and measure operators by using the basis expansion of the infinite variables (i.e., with operators applied to the basis functions).  In certain cases, the expansion coefficients can be useful in evaluating certain measure types; for example, the coefficients will correspond to the statistical moments of the infinite variables when PCE is applied to a stochastic formulation with a basis that is orthogonal to the probability density function and these moments are often used to express expectations and risk measures \cite{zymler2013distributionally}.

After basis expansion representations are incorporated, the problem is fully defined in terms of the finite-dimensional variables $\tilde{y}_i$ and $z$. However, this representation is not yet tractable, since it still contains infinite-dimensional objects (e.g., basis functions and associated operators).  To deal with this, we consider the residual (i.e., finite error) $R(d)$ associated with performing this projection on each constraint and on the objective. Each resulting residual will be made as small as possible by exacting that they be orthogonal to a set of weight functions $\psi_k(d),\; k \in \mathcal{K}, d\in \mathcal{D}$:
\begin{equation}
    \langle R, \psi_k \rangle_w = 0, \ \forall k \in \mathcal{K}
    \label{eq:mwr_residual}
\end{equation}
where $\langle \cdot,\cdot \rangle_w$ denotes the inner product between functions using the appropriate weighting function $w(d)$ for the given space:
\begin{equation}
    \langle R, \psi_k \rangle_w = \int_{d' \in \mathcal{D}}R(d')\psi_k(d')w(d')dd'.
    \label{eq:inner_product}
\end{equation}
The weight functions are typically chosen such that $|\mathcal{K}| = |\mathcal{I}|$. The projection results in a tractable finite-dimensional formulation; 
what remains is our choice of the weight functions $\psi_k(\cdot)$. This choice gives rise to a range of techniques; if the Galerkin method is applied then we choose $\phi_k(d) = \psi_k(d)$ and have that $\mathcal{I} = \mathcal{K}$. This induces the first $|\mathcal{I}|$ terms of the residuals in the trial functions to vanish if the functions are orthogonal \cite{graham2013modeling}. Another popular choice is that of orthogonal collocation, where we choose $\psi_k(d) = \delta(d - \hat{d}_k)$;  here, the set $\hat{d}_k,\; k \in \mathcal{K}$ denote collocation points (i.e., particular infinite parameter supports) and $\delta(\cdot)$ is the Dirac delta function. This approach seeks to enforce that the residual is zero at the collocation points. When orthogonal basis functions are chosen and this is applied over a set of transcription points (i.e., finite elements), we obtain a method known as orthogonal collocation over finite elements. A variety of other methods such as least squares and the method of moments can also be employed to weight the residuals (these are discussed in detail in \cite{finlayson2013method}).

The transformation of \eqref{eq:infinite_form} to a finite-dimensional form via MWR is, in general, a nonlinear transformation (depending on the choices of the trial functions $\phi_i(\cdot)$, weight functions $\psi_k (\cdot)$, and their corresponding space). However, there exist special cases where the transformation is linear, as is often the case with PCE transformations \cite{muhlpfordt2019chance}. Advantages of employing MWR instead of direct transcription is that one obtains functional representations for the infinite variables (as opposed to values at the support points), one can achieve better stability for boundary-valued problems, and one can obtain better accuracy for certain formulations \cite{devolder2010solving}. On the other hand, the main disadvantage of MWR is that evaluating differential and measure operators and inner products tends to be cumbersome (especially for nonlinear formulations). Also, basis functions can be difficult to derive for formulations with multivariate infinite domains. In our abstraction, we provide the modeling elements that facilitate the implementation of these transformation techniques.

\subsection{Transformation Framework in \texttt{InfiniteOpt.jl}} \label{sec:infiniteopt_solution}

In Section \ref{sec:infiniteopt_modeling} we discussed how our unifying abstraction is implemented in \texttt{InfiniteOpt.jl}; one creates a model as an   \texttt{InfiniteModel} object. In this section, we discuss a general transformation framework incorporated into \texttt{InfiniteOpt.jl} that facilitates the implementation of different transformation approaches (e.g., direct transcription and MWRs). We also outline the efficient direct transcription capabilities that are currently implemented. 

\begin{figure}[!htb]
    \centering
    \includegraphics[width=0.6\textwidth]{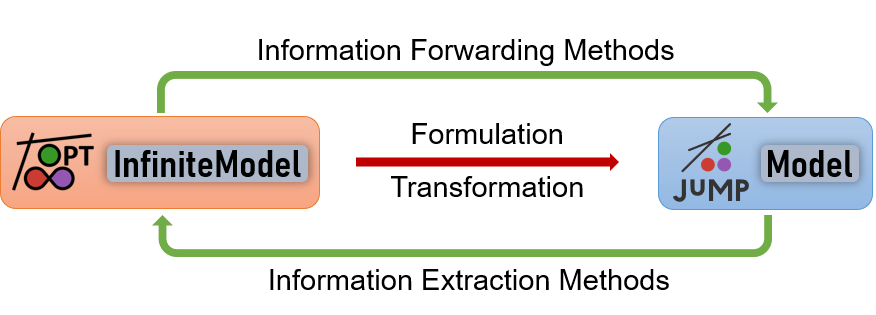}
    \caption{Transformation framework employed by \texttt{InfiniteOpt.jl} for converting an \texttt{InfiniteModel} into a  \textrm{JuMP.jl} \texttt{Model}.}
    \label{fig:optimizer_model}
\end{figure}

The framework centers around applying a transformation to the \texttt{InfiniteModel} that converts it to a standard \textrm{JuMP.jl} \texttt{Model} object (referred to as an optimizer model in this context). The optimizer model can then be solved using the optimizers implemented in \textrm{MathOptInterface.jl} \cite{legat2020mathoptinterface}. Moreover, this framework features a collection of methods to enable a seamless interface between the \texttt{InfiniteModel} and its corresponding optimizer model to facilitate capabilities such as information extraction (e.g., solution queries) that do not require direct interrogation of the optimizer model. This framework is summarized in Figure \ref{fig:optimizer_model}. This software structure distinguishes \texttt{InfiniteOpt.jl} from other software tools (e.g., \texttt{Pyomo.dae} and \texttt{Gekko}) whose implementations center around (and are limited to) direct transcription. Thus, \texttt{InfiniteOpt.jl} is unique in providing a flexible API for solving InfiniteOpt problems. 

Following this framework, a desired solution scheme is incorporated by defining a few prescribed methods (principally \texttt{build\_optimizer\_model!}) to implement the associated transformation. This methodology is implicitly invoked on line \ref{line:optimize} of Code Snippet \ref{code:infiniteopt_example} where \texttt{optimize!} creates an optimizer model using the prescribed transformation and then solves it with the desired optimizer. The full technical detail of this API is beyond the scope of this work and is available via the \texttt{InfiniteOpt.jl} documentation.

\begin{figure}[!htb]
    \centering
    \includegraphics[width=0.6\textwidth]{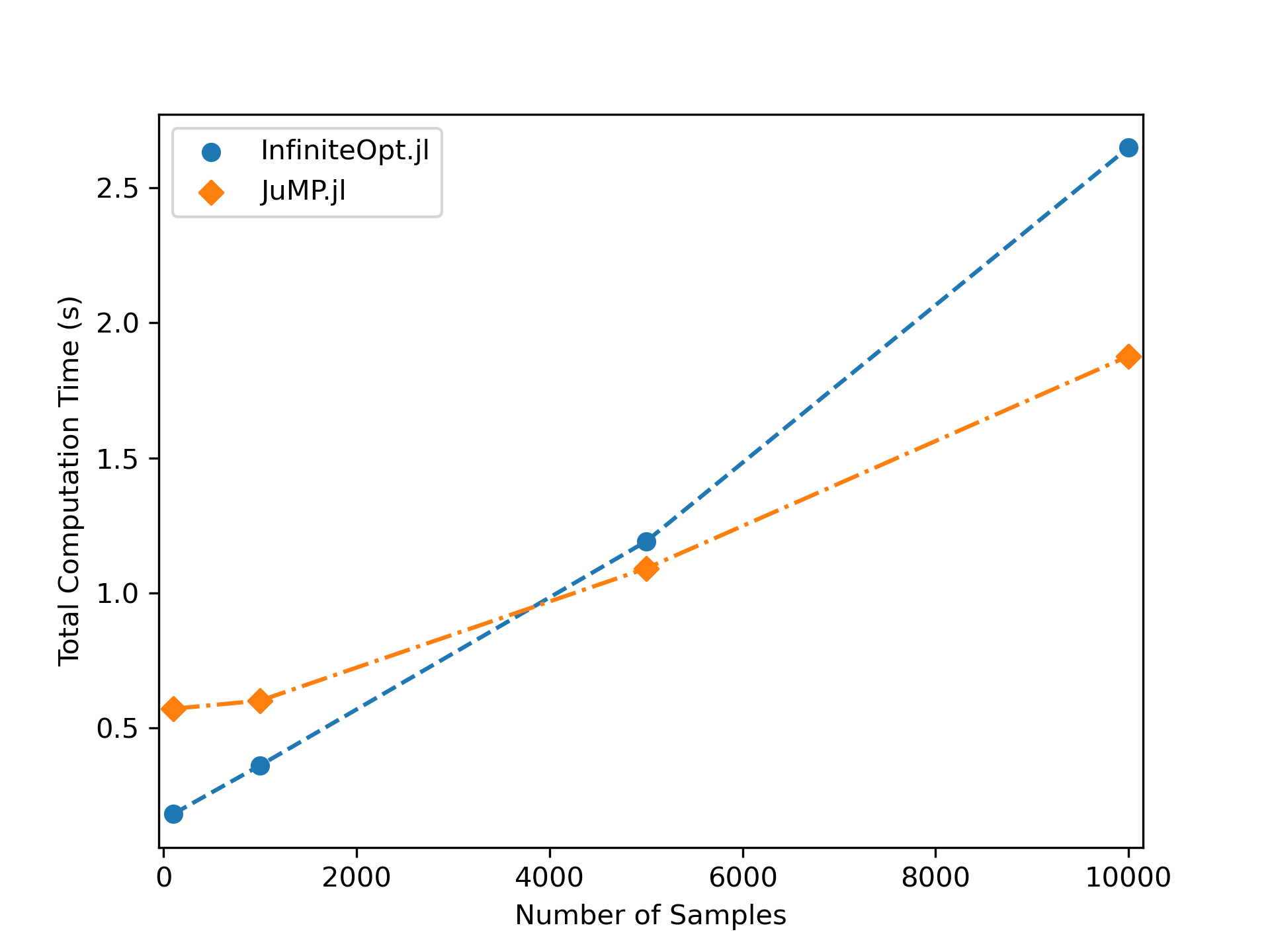}
    \caption{Juxtaposition of the total computation time used to formulate and solve a stochastic optimization problem \cite{pulsipher2019scalable} using MC sampling implemented in \textrm{InfiniteOpt.jl v0.4.1} and manual transcription in \textrm{JuMP.jl v0.21.8}.}
    \label{fig:transcriptionopt_scale}
\end{figure}

\texttt{InfiniteOpt.jl} provides an efficient implementation of direct transcription following the procedures described in Section \ref{sec:transcription}; this serves as the default transformation technique for InfiniteOpt models. These techniques are implemented in a sub-module called \textrm{TranscriptionOpt} that follows the optimizer model framework shown in Figure \ref{fig:optimizer_model}. The \textrm{TranscriptionOpt} module features a sophisticated finite support generation and management system that enables tackling a wide variety of infinite-dimensional optimization formulations using diverse evaluation techniques for measure and derivative operators. Moreover, its automatic transcription is efficient and compares competitively to manually transcribing a problem and implementing it via \texttt{JuMP.jl}. This incredible behavior is demonstrated in Figure \ref{fig:transcriptionopt_scale} where a two-stage stochastic optimization problem (the 3-node distribution network example featured in \cite{pulsipher2019scalable}) is solved for a range of MC samples using automatic transcription in \texttt{InfiniteOpt.jl} and manual transcription in \texttt{JuMP.jl}. We note that, contrary to other software implementations, automatic transcription in \texttt{InfiniteOpt.jl} denotes a minor computational expense relative to manual transcription with the benefit of avoiding the errors commonly incurred by transcribing InfiniteOpt formulations manually. 

\section{Innovations Enabled by Unifying Abstraction} \label{sec:consequences}

In this section, we discuss innovative modeling approaches that are enabled by the proposed InfiniteOpt abstraction. Specifically, these innovations are facilitated by the ability to transfer modeling techniques and constructs across disciplines. For instance, in the proposed abstraction, there is no explicit notion of time, space, or random domains (all domains are treated mathematically in the same way); as such, one can easily identify analogues of modeling elements across disciplines. 

\subsection{Measure Operators} \label{sec:generalized_measures}

Here we highlight the advantages of using an abstraction that focuses on measure operators. We will highlight how InfiniteOpt formulations from different disciplines are connected via analogous mathematical features; we place particular attention to common features arising in dynamic and stochastic optimization problems.  

\subsubsection{Expectation Measures} \label{sec:expect}

\begin{figure}[!htb]
    \centering
    \includegraphics[width=0.5\textwidth]{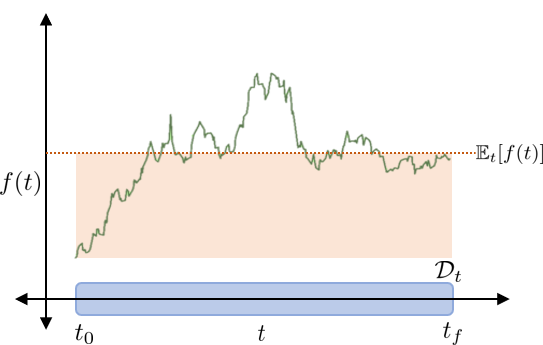}
    \caption{Visualization of the expectation measure $\mathbb{E}_t[f(t)] = \frac{1}{t_f - t_0}\int_{t_0}^{t_f} f(t) dt$ where the rectangle formed has an area equal to that of the region under $f(t)$.}
    \label{fig:time_expect}
\end{figure}

In Section \ref{sec:measures}, we provide examples of typical measures used to formulate objective functions in InfiniteOpt problems (space-time integrals and expectations). For simplicity in the presentation, we consider a temporal domain $\mathcal{D}_t$ and consider the time integral: 
\begin{equation}
    \int_{t \in \mathcal{D}_t} f(t) dt
    \label{eq:optimal_control_1d}
\end{equation}
where we write $f(t)=f(y(t), t)$ for compactness. Minimizing the measure \eqref{eq:optimal_control_1d} seeks to shape the cost function (a surface defined in the domain $\mathcal{D}_t$) in a way that it displaces the entire surface. Analogously, minimizing the expectation measure:
\begin{equation}
  \mathbb{E}_\xi[f(\xi)]:=  \int_{\xi \in \mathcal{D}_\xi} f(\xi) p(\xi)d\xi
    \label{eq:expectation_1d}
\end{equation}
 shapes the cost surface (defined in the domain $\mathcal{D}_\xi$) in a way that displaces the surface. An obvious difference between time integral and the random expectation is that the random expectation is weighted by a probability density function $p(\cdot)$ and this gives flexibility to put more/less emphasis on different locations of the random domain $\mathcal{D}_\xi$. Thus, as a simple example on how one can transfer modeling strategies, we can formulate the expectation over the temporal domain using a weighting function as:
\begin{equation}
    \mathbb{E}_t[f(t)] := \int_{t \in \mathcal{D}_t} f(t) w(t) dt.
    \label{eq:optimal_control_weighted}
\end{equation}
Note that the selection of the notation $t$ to denote the infinite domain is arbitrary; one can simply define a general infinite parameter $d$. If one defines $w(t)=1/S$ with $S=\int_{t\in\mathcal{D}_t}dt$, we can see that the above measure is the time average of $f(t)$ (with equal weight placed at each location in the domain). If the time domain is $D_{t}=[t_0,t_f]$, we have that $S=t_f-t_0$. Figure \ref{fig:time_expect} provides a geometric interpretation of the time-expectation; here, the area of the rectangle with height $\mathbb{E}_t[f(t)]$ and width $S$ is equivalent to the area under $f(t)$ \cite{stewart2009calculus}. This is just a scaled version of the integral measure \eqref{eq:optimal_control_1d} and will thus shape $f(t)$ in similar manner. In fact, the integral \eqref{eq:optimal_control_1d} is just a special case of measure \eqref{eq:optimal_control_weighted} (obtained by setting $w(t) = 1$); also, note that this approach is equivalent to using a weighting function $w(t)$ that corresponds to the probability density function of a uniform random parameter $\xi \sim \mathcal{U}(t_0, t_f)$. As such, one can envision defining weighting functions associated with probability densities of different random parameters (e.g., Gaussian, exponential, Weibull); this would have the effect of inducing interesting prioritization strategies that can be used to shape the cost surface in desirable ways. For instance, a Gaussian weighting function places emphasis on the middle of the time domain (and emphasis decays rapidly as one moves away from the center of the domain), while an exponential weighting function places emphasis at the beginning of the domain (and decays rapidly as one marches in time).  This modeling feature can be useful in dynamic optimization and optimal control problems in which it is often desirable to place more/less emphasis on initial or final conditions.  For instance, in infinite-horizon problems, $w(\cdot)$ obtained from the probability density function of an exponential density function behaves as a discount factor \cite{petrik2008biasing,shin2021exponential}. 

\subsubsection{Risk Measures} \label{sec:risk_measures}

A large collection of risk measures have been proposed in the stochastic optimization literature to manipulate random cost functions in desirable ways (e.g., to minimize impacts of extreme events)  \cite{ruszczynski2006optimization}. In this context, risk measures are typically interpreted as measures that aim to shape the tail of the probability densities of cost or constraint functions.  A popular risk measure used for this purpose is the conditional value-at-risk (CVaR):
\begin{equation}
    \text{CVaR}_\xi(f(\xi); \alpha) := \min_{\hat{f} \in \mathbb{R}} \left\{\hat{f} + \frac{1}{1 - \alpha} \mathbb{E}_\xi\big(\max(f(\xi) - \hat{f}, 0)\big) \right\}
    \label{eq:cvar_def}
\end{equation}
where $\alpha \in [0,1)$ is a desired probability level and $\hat{f}$ is an auxiliary variable.
\\

One can show that the value of the auxiliary variable that minimizes the inner function of CVaR is given by $\hat{f}^*=Q_\xi(f(\xi); \alpha)$, which is the $\alpha$-quantile of $f(\xi)$ \cite{rockafellar2000optimization}.  We recall that the quantile is defined as:
\begin{equation}
    Q_\xi(f(\xi); \alpha):=\inf_{\hat{f} \in \mathbb{R}} \left\{\mathbb{P}_\xi\left(f(\xi) \leq \hat{f}\right) \geq \alpha \right\}.
\end{equation}
The quantile is thus the threshold value for $f(\xi)$ such that the probability of finding this function below the threshold is at least $\alpha$.
\\

One can also show that  CVaR is a conditional expectation of the form: 
\begin{equation}
    \text{CVaR}_\xi(f(\xi); \alpha) = \mathbb{E}_\xi\big[f(\xi) \ | \ f(\xi) \geq Q_\xi(f(\xi); \alpha)\big]
    \label{eq:cvar_easy}
\end{equation}
Hence, minimizing CVaR has the effect of minimizing the conditional expectation over the $1-\alpha$ probability region with the highest cost, thus hedging against extreme events. Moreover, the calculation of CVaR at a probability level $\alpha$ implicitly defines the quantile of $f(\xi)$.

\begin{figure}[!htb]
    \centering
    \includegraphics[width=0.5\textwidth]{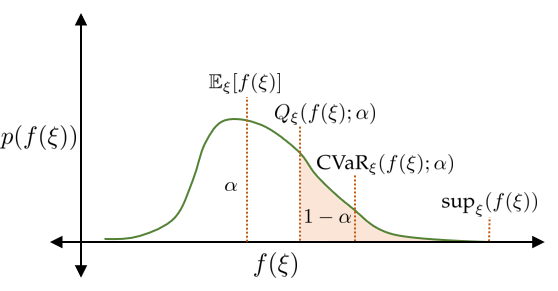}
    \caption{An illustration of $\text{CVaR}_\xi(f(\xi); \alpha)$ in terms of the probability density function $p(f(\xi))$.}
    \label{fig:cvar1}
\end{figure}

A key property of CVaR is that it is a general measure that captures the expectation $\text{CVaR}_\xi(f(\xi); \alpha) = \mathbb{E}_\xi[f(\xi)]$ as $\alpha \rightarrow 0$ and the worst-case $\text{CVaR}_\xi(f(\xi); \alpha) = \sup_\xi(f(\xi))$ as $\alpha \rightarrow 1$. Figure \ref{fig:cvar1} shows how these measures are typically interpreted in terms of the probability density function of the cost $f(\xi)$, motivating the $\alpha$ and $1-\alpha$ probability regions.

\begin{figure}[!htb]
    \centering
    \includegraphics[width=0.4\textwidth]{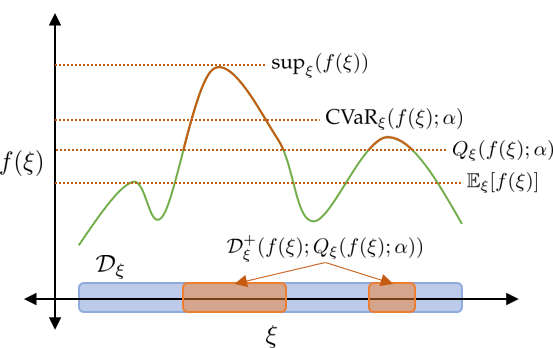}
    \caption{Illustration of $\text{CVaR}_\xi(f(\xi); \alpha)$ following the representation given in \eqref{eq:cvar_expect}. This provides an alternative view of the probabilistic representation shown in Figure \ref{fig:cvar1}.}
    \label{fig:cvar2}
\end{figure}

Through the perspective of the proposed unifying abstraction, one can interpret CVaR as a measure that captures the {\em excursion} of a function (field) from a given threshold. To see this, we define the positive and negative function excursion sets of $f(\xi)$ which denote the range of $f(\xi)$ above and below a threshold $\hat{f}$, respectively:
\begin{equation}
    \begin{aligned}
        &\mathcal{D}^+_\xi(f(\xi);\hat{f}) := \{\xi \in \mathcal{D}_\xi: \ f(\xi) \geq \hat{f}\} \\
        & \mathcal{D}^-_\xi(f(\xi); \hat{f}) := \{\xi \in \mathcal{D}_\xi: \ f(\xi) \leq \hat{f}\}
    \end{aligned}
    \label{eq:function_excursions}
\end{equation}
where $\mathcal{D}^+_\xi(f(\xi);\hat{f}) \subseteq \mathcal{D}_\xi$ and $\mathcal{D}^-_\xi(f(\xi);\hat{f}) \subseteq \mathcal{D}_\xi$ are the positive and negative function excursion sets, respectively. We simplify the notation for these excursion sets by using $\mathcal{D}^-_\xi(\hat{f})$ and $\mathcal{D}^+_\xi(\hat{f})$ . Using these definitions, we can express the quantile function $q_\alpha:=Q_\xi(f(\xi); \alpha)$ as:
\begin{equation}
    q_\alpha= \inf_{\hat{f} \in \mathbb{R}} \left\{\int_{\xi \in \mathcal{D}^-_\xi(\hat{f})} p(\xi) d\xi \geq \alpha \right\}.
    \label{eq:gen_quantile}
\end{equation}
This reveals that CVaR considers the expectation of $f(\xi)$ over the restricted domain $\mathcal{D}^+_\xi(q_\alpha)$, which indexes the $1-\alpha$ probabilistic region shown in Figure \ref{fig:cvar1}:
\begin{equation}
    \text{CVaR}_\xi(f(\xi); \alpha) = \frac{1}{1-\alpha} \int_{\xi \in  \mathcal{D}^+_\xi(q_\alpha)} f(\xi) p(\xi) d\xi.
    \label{eq:cvar_expect}
\end{equation}
Figure \ref{fig:cvar2} illustrates that $\text{CVaR}_\xi(f(\xi); \alpha)$ using this functional interpretation for a realization of $f(\xi)$ and compares it to the other measures shown in Figure \ref{fig:cvar1}.

\begin{figure}[!htb]
    \centering
    \includegraphics[width=0.4\textwidth]{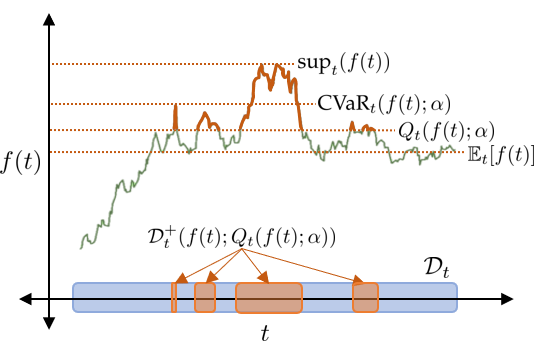}
    \caption{Illustration of $\text{CVaR}_t(f(t); \alpha)$, as represented in \eqref{eq:cvar_time}.}
    \label{fig:cvar3}
\end{figure}

The representation of CVaR shown in \eqref{eq:cvar_expect} highlights that this can be seen as a measure operator that can be applied to different types of infinite domains (e.g., space and time). For instance, shaping time-dependent trajectories to minimize extreme values is applicable to dynamic optimization (e.g., to prevent peak costs or extreme operating conditions that can hinder safety). For instance, Risbeck and Rawlings \cite{risbeck2019economic} recently proposed a model predictive control formulation that has the dual goal of minimizing total accumulated cost and peak cost:
\begin{equation}
    \int_{t \in \mathcal{D}_t}f(t)dt + \lambda \cdot \mathop{\text{max}}_{t\in \mathcal{D}_t}\; f(t)
    \label{eq:peak_cost}
\end{equation}
where $\lambda \in \mathbb{R}_+$ is a trade-off parameter. More generally, we might want to penalize a subdomain of the time domain that attains the highest costs; this can be achieved by applying a CVaR measure operator to the time-dependent cost $f(t)$: 
\begin{equation}
    \text{CVaR}_t(f(t); \alpha) = \frac{1}{1-\alpha} \int_{t \in  \mathcal{D}^+_t(q_\alpha))} f(t) p(t) dt
    \label{eq:cvar_time}
\end{equation}
where the density function can be selected as $p(t) = \frac{1}{t_f - t_0}$. Note that this definition implicitly introduces the notion of a quantile in a time domain; this quantile is given by:
\begin{equation}
    q_\alpha:=Q_t(f(t);\alpha) = \inf_{\hat{f} \in \mathbb{R}} \left\{\int_{t \in \mathcal{D}^-_t(f(t); \hat{f})} p(t) dt \geq \alpha \right\}.
    \label{eq:time_quantile}
\end{equation}

Using the properties of CVaR we can conclude that: 
\begin{equation}
    \begin{aligned}
        &\lim_{\alpha \rightarrow 0} \text{CVaR}_{t}(f(t); \alpha) = \int_{t \in \mathcal{D}_t} f(t)p(t) dt  \\
        &\lim_{\alpha \rightarrow 1} \text{CVaR}_{t}(f(t); \alpha) = \mathop{\text{max}}_{t\in \mathcal{D}_t}\; f(t).
    \end{aligned}
\end{equation}
This highlights that CVaR provides an intuitive measure that can help shape time-dependent trajectories. Figure \ref{fig:cvar3} illustrates the application of this measure over the time domain and shows that this is analogous to the application over a random domain shown in Figure \ref{fig:cvar2}. We highlight that, in practice, $\text{CVaR}_t(f(t); \alpha)$ is computed by using \eqref{eq:cvar_def} (defined for $t$ instead of $\xi$). 
\\

This example illustrates how the CVaR construct used in stochastic optimization can be utilized to motivate new formulations for other optimization disciplines. The measure-centric unifying abstraction facilitates this process by capturing objects via measure operators. The process of transferring ideas through the measure-centric abstraction is also amendable to other risk measures (see \cite{krokhmal2013modeling} for a review of risk measures). We present a numerical study of using CVaR for an optimal control formulation in Section \ref{sec:cases_cvar}. 

\subsubsection{Event Constraints} \label{sec:event_constrs}

An interesting modeling paradigm that arises from the proposed abstraction are {\em event constraints}. These constraints generalize the notion of chance constraints, excursion set conditions, and exceedance probabilities that are used in different scientific disciplines (e.g., stochastic optimization, reliability engineering, and random fields). This unified view also promotes transfer of modeling concepts across disciplines; for instance, we will see that chance constraints in a random domain are analogous to exceedance times in a time domain. 
\\

To exemplify the notion of event constraints, consider the so-called {\em excursion time}; this measure is widely used in reliability analysis of dynamical systems and is defined as the fraction of time that a function $h(t),\; t\in \mathcal{D}_t$ is above a given threshold \cite{au2001first}. Here, we consider a zero threshold value to give the event constraint: 
\begin{equation} 
    \mathbb{P}_t\left( \{\exists t \in \mathcal{D}_t : h(t) > 0\}\right) \leq \alpha.
    \label{eq:time_excursion}
\end{equation}
where $\alpha\in [0,1]$. In the context of an InfiniteOpt problem, the function $h(t),\,t\in \mathcal{D}_t$ can denote a constraint $h(y(t), z, t),\,t\in \mathcal{D}_t$. The excursion time is expressed as a probability-like measure of the form: 
\begin{equation} 
    \mathbb{P}_t\left(\{\exists t \in \mathcal{D}_t : h(t) > 0\}\right)=\int_{t\in \mathcal{D}_t} \mathbbm{1}\left[\{\exists t \in \mathcal{D}_t : h(t) > 0\}\right]w(t)dt
    \label{eq:time_excursion2}
\end{equation}
where $w:\mathcal{D}_t\to [0,1]$ is a weighting function satisfying $\int_{t\in \mathcal{D}_t}w(t)dt=1$. The excursion time measure can be interpreted as the {\em fraction of time} under which the event of interest occurs;  for instance, in safety analysis, one might be interested in determining the fraction of time that a constraint is violated and to ensure that this is not greater than some fraction $\alpha$.  Alternatively, we could also search to minimize this measure (by using it as an objective). The excursion time constraint is an event constraint that can help shape a time-dependent trajectory in interesting and non-intuitive ways.
\\

One can generalize the excursion time measure by constructing complex events. For instance, consider that we want to determine the fraction of time that any of the time-varying constraints $h_k(t),\; k\in \mathcal{K}$ crosses a threshold. This can be expressed as:
\begin{equation} 
    \mathbb{P}_t\left(\bigcup_{jk\in \mathcal{K}} \{\exists t \in \mathcal{D}_t : h_k(t) > 0\}\right).
    \label{eq:time_excursion3}
\end{equation}
Here, $\bigcup$ is a logical \texttt{or} operator. If we want to determine fraction of time that all constraints are violated then we use: 
\begin{equation} 
    \mathbb{P}_t\left(\bigcap_{h \in \mathcal{K}} \{\exists t \in \mathcal{D}_t : h_k(t) > 0\}\right).
    \label{eq:time_excursion4}
\end{equation}
Here, $\bigcap$ is a logical \texttt{and} operator. 
\\

\begin{figure}[!htb]
    \centering
    \begin{subfigure}[b]{0.3\textwidth}
        \centering
        \includegraphics[width=\textwidth]{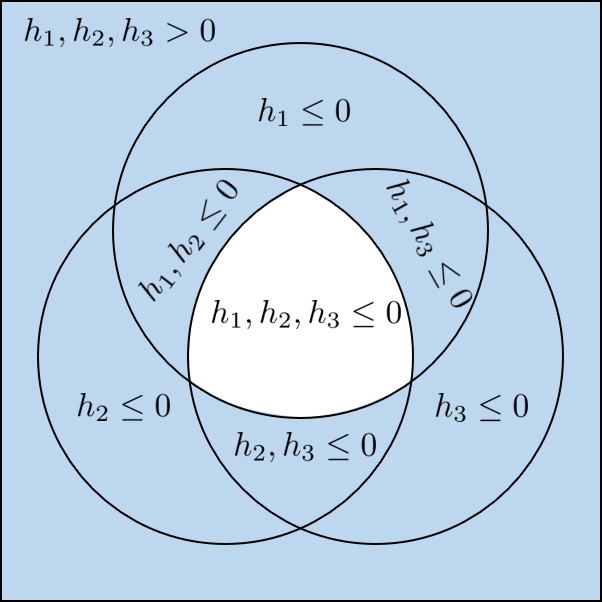}
        \caption{Excursion Constraint}
        \label{fig:excursion_logic}
    \end{subfigure}
    \quad \quad \quad
    \begin{subfigure}[b]{0.3\textwidth}
        \centering
        \includegraphics[width=\textwidth]{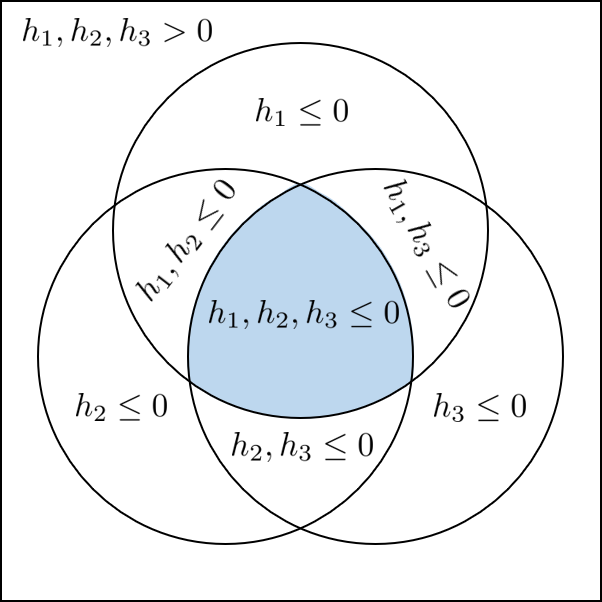}
        \caption{Joint-Chance Constraint}
        \label{fig:joint_logic}
    \end{subfigure}
    \caption{Logical event regions (shown in blue) constrained by the event constraints \eqref{eq:chance_excursion} and \eqref{eq:joint_chance5}. In particular, they constrain the condition $h_1 >0 \cup h_2 > 0 \cup h_3 > 0$ and the condition $h_1 \leq 0 \cap h_2 \leq 0 \cap h_3 \leq 0$, respectively.}
    \label{fig:logic_regions}
\end{figure}

To highlight transfer across different disciplines, we recognize that the excursion time is directly analogous to a chance constraint (operating in the random space, as opposed to time) and our previous analysis suggests that one can construct chance constraints that capture complex events. For instance, consider the event constraint: 
\begin{equation} 
    \mathbb{P}_\xi\left(\bigcup_{k \in \mathcal{K}} \{\exists \xi \in \mathcal{D}_\xi : h_k(\xi) > 0\}\right) \leq \alpha.
    \label{eq:chance_excursion}
\end{equation}
Here, we see that the constraint is directly analogous to the event constraint \eqref{eq:time_excursion3}. Figure \ref{fig:excursion_logic} illustrates the logical event space that constraint \eqref{eq:chance_excursion} shapes. We thus see that events can be defined in a general form over different infinite domains; to highlight this fact, we consider the event constraint:
\begin{equation} 
    \mathbb{P}_d\left(\bigcup_{k \in \mathcal{K}} \{\exists d \in \mathcal{D} : h_k(d) > 0\}\right) \leq \alpha. 
    \label{eq:chance_excursiond}
\end{equation}

The previous event constraint is different to traditional joint-chance constraints used in stochastic optimization: 
\begin{equation}                             
    \mathbb{P}_\xi\left(\bigcap_{k \in \mathcal{K}}\{\exists \xi \in \mathcal{D}_\xi : h_k(\xi) \leq 0 \}\right) \leq \alpha.
    \label{eq:joint_chance5}
\end{equation}
This makes it more readily apparent that traditional joint-chance constraints consider the logical event space that is complementary to that of constraint \eqref{eq:chance_excursion}. Figure \ref{fig:joint_logic} shows this region which is the logical complement of Figure \ref{fig:excursion_logic}. 
\vspace{0.2in}

\begin{figure}[!htb]
    \centering
    \includegraphics[width=0.3\textwidth]{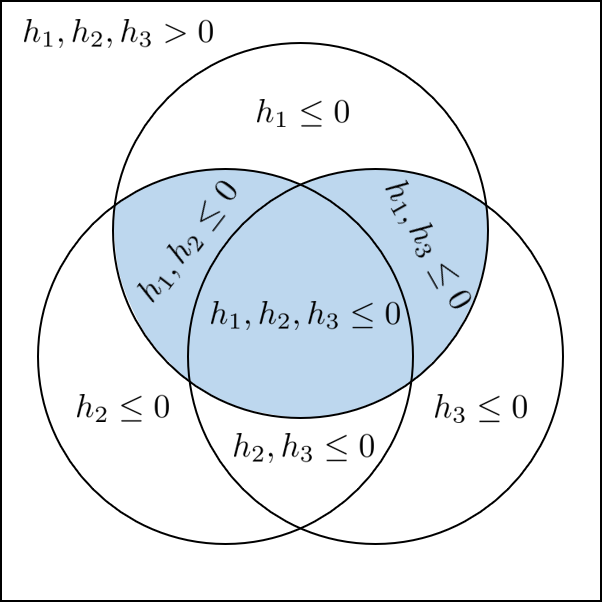}
    \caption{Illustration of logical event region captured by constraint \eqref{eq:novel_chance}.}
    \label{fig:novel_logic}
\end{figure}

Another interesting insight that we draw from event constraints is that logical operators (e.g., $\cap$ and $\cup$) can be used to model complex decision-making logic. For example, following the constraint system shown in Figure \ref{fig:logic_regions}; we might consider the logical event region derived from the condition that $h_1 \leq 0 \cap (h_2 \leq 0 \cup h_3 \leq 0)$ giving the event constraint:
\begin{equation}                             
    \mathbb{P}_\xi\left( \{\forall \xi \in \mathcal{D}_\xi : h_1 \leq 0 \cap (h_2 \leq 0 \cup h_3 \leq 0)\}\right) \geq \alpha.
    \label{eq:novel_chance}
\end{equation}
This is depicted in Figure \ref{fig:novel_logic}; we note that this event constraint encapsulates a wider probabilistic region relative to that of the traditional joint-chance constraint \eqref{eq:joint_chance5}.
\\

In summary, the presented examples illustrate that excursion time constraints and chance constraints as special cases of event constraints. This crossover also led us to representing joint-chance constraints with logical operators which introduce the notion of embedding problem-specific logic to shape the probabilistic region these constraints consider. We illustrate this example further with a stochastic optimal power flow case study in Section \ref{sec:sopf}. 

\subsection{Random Fields} \label{sec:random_fields}

\begin{figure}[!htb]
    \centering
    \includegraphics[width=0.4\textwidth]{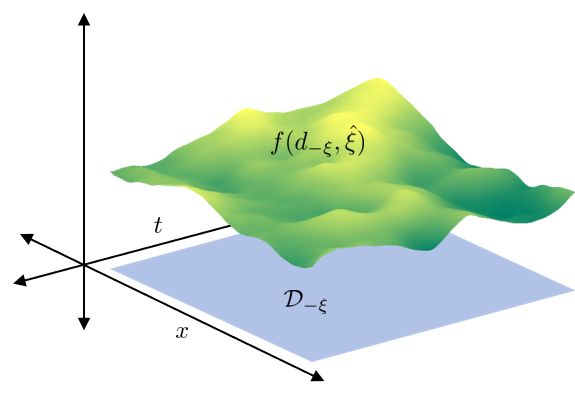}
    \caption{Sample of a random field defined over space and time.}
    \label{fig:random_field}
\end{figure}

We now discuss how the abstraction inspires the incorporation of modeling concepts from random field theory into InfiniteOpt formulations.  A random field is a random function with realizations of the form $f(d) : \mathcal{D} \mapsto \mathbb{R}^{n_f}$ \cite{adler2010geometry}. For instance, one can think of a continuous-time trajectory that is random (e.g., due to uncertainty in a differential equation model that generates it). A random field generalizes the notion of multivariate random variables (e.g., a multivariate Gaussian $\xi \sim\mathcal{N}(\mu, \Sigma)$) that are jointly-distributed to that of infinite jointly-distributed random variables that are indexed over some continuous domain. Another example of a random field is that of a dynamic Gaussian process $\xi(t) \sim \mathcal{GP}(\mu(t), \Sigma(t, t'))$ for $t\in \mathcal{D}_t$. Figure \ref{fig:random_field} depicts a realization of a random field. 

Modeling concepts from random field theory allows us to incorporate dependencies of a random parameter over other infinite domains (e.g., account for stochastic processes). For instance, in the context of our abstraction, consider the random time-dependent function $y(t, \xi)$. This decision function is random due to the presence of $\xi \in \mathcal{D}_\xi$ and is also indexed over $t \in \mathcal{D}_t$; this means that it can be precisely characterized as a random field. From this it follows that all infinite variables $y(d_{-\xi}, \xi)$ are random fields. Hence, optimally choosing $y(d_{-\xi}, \xi)$ in a given formulation amounts to engineering a random field that is in general correlated over the infinite domain $\mathcal{D}_{-\xi}$. This important observation enables us to connect a wide breadth of optimization disciplines to random field theory. For instance, the theory and methods behind random field constructs like excursion sets, excursion probabilities, and expected Euler characteristics are amendable to the measure operators and/or event constraints formalized in our abstraction (see Sections \ref{sec:risk_measures} and \ref{sec:event_constrs}) \cite{adler2000excursion}.

\begin{figure}[!htb]
    \centering
    \includegraphics[width=0.4\textwidth]{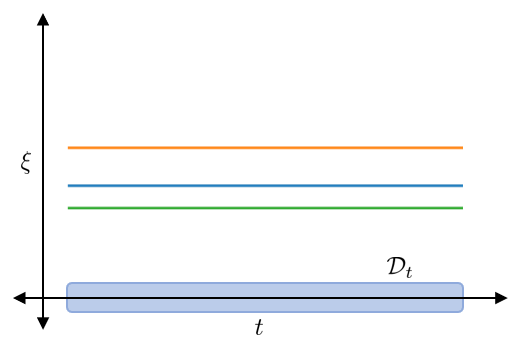}
    \caption{Realizations of a random parameter $\xi \in \mathcal{D}_\xi$ that are invariant relative to other infinite domains (e.g., time $t$).}
    \label{fig:static_uncertainty}
\end{figure}

\begin{figure}[!htb]
    \centering
    \begin{subfigure}[b]{0.4\textwidth}
        \centering
        \includegraphics[width=\textwidth]{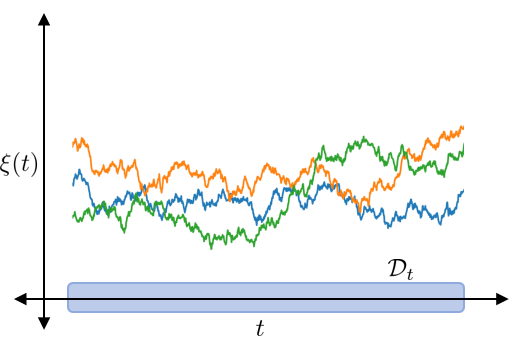}
        \caption{Correlated}
        \label{fig:correlated_field}
    \end{subfigure}
    \quad \quad \quad
    \begin{subfigure}[b]{0.4\textwidth}
        \centering
        \includegraphics[width=\textwidth]{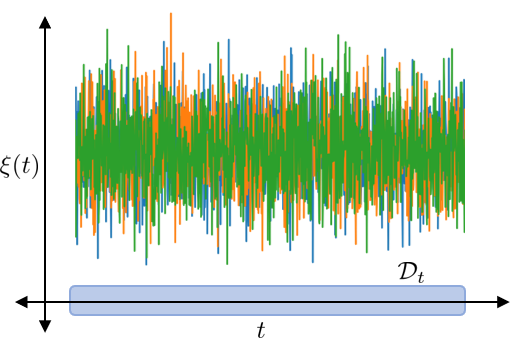}
        \caption{Uncorrelated}
        \label{fig:uncorrelated_field}
    \end{subfigure}
    \caption{Realizations of a time-dependent random field $\xi(t)$ with and without time correlation.}
    \label{fig:field_uncertainty}
\end{figure}

The particular class of infinite variables discussed up to this point only consider a random parameter $\xi$ that is invariant over other domains, as exemplified for a few realizations of $\xi$ in Figure \ref{fig:static_uncertainty}. This means that, although the infinite variables in this class of formulations are general random fields, the input uncertainty $\xi$ is taken to be {\em static} random parameter that does not capture any spatial or temporal correlation. We can overcome this modeling limitation by extending our abstraction to consider infinite parameter functions (e.g., random parameters that are functions of other infinite parameters). In the context of our current example, this is accomplished by defining the infinite parameter function $\xi(t) \in \mathcal{D}_{\xi(t)},\; t\in \mathcal{D}_t$ which denotes a known random field whose sample domain $\mathcal{D}_{\xi(t)}$ is a set of temporal functions. With this extension, we can now define the infinite variable $y(t, \xi(t))$, which denotes a random field variable (an infinite variable) where $\xi(t)$ is a known random field that can capture correlations of $\xi$ over time $t$. Note that $\xi(t)$ is explicitly shown as an input to $y$ so that it is distinguished from deterministic infinite-dimensional variables $y(t)$. Figure \ref{fig:field_uncertainty} shows a realization of $\xi(t)$ in the case that it exhibits temporal correlation and in the case that $\xi$ approaches no temporal correlation. Note this can also capture the static case shown in Figure \ref{fig:static_uncertainty}. 

Random fields have been widely used in disciplines such as robust topology design optimization \cite{zhao2014robust, zhang2018robust}; however, to the best of our knowledge, the integration of random field theory in optimization represents a new class of optimization formulations. We call this class {\em random field optimization} and we note that this is a special case of the proposed InfiniteOpt abstraction. This class provides effective means of capturing correlation of uncertainty over other domains an optimization formulation which is enabled by the rich theory that has been developed for random fields \cite{adler2010geometry}. Moreover, random field optimization can be seen as a {\em generalization of stochastic optimization}. For instance, it provides an intriguing alternative to multi-stage (discrete time) stochastic formulations which cannot be readily generalized to continuous time or general infinite domains (e.g., space) \cite{shapiro2003inference}, whereas our proposed formulation class is defined for general infinite domains (e.g., space and time) and can incorporate random field uncertainty models such as Gaussian processes.
\\

Another interesting observation is that extending the proposed abstraction to include infinite parameter functions can also be done for non-random infinite parameters. For instance, we could account for a time-dependent spatial parameter $x(t) \in \mathcal{D}_{x(t)}$ where $\mathcal{D}_{x(t)}$ is some space. We leave a rigorous formalization of such an extension to our abstraction and its implications for random field optimization problems to future work. 

\subsection{Problem Analysis} \label{sec:characterizations}

Here we highlight some of the advantages that arise from characterizing InfiniteOpt problems directly in accordance with formulation \eqref{eq:infinite_form}, in contrast to the standard practice of expressing them solely via finite reformulations. For instance, within the area of dynamic optimization, it is commonplace to abstract and formalize problem classes in discrete time (i.e., in a transcribed form) \cite{rawlings2017model}. This practice tends to make problems more difficult to formulate since they are inherently coupled with the transformation scheme employed (e.g., orthogonal collocation or explicit Euler). Hence, decoupling the problem definition from its transformation helps to ease its formalization. Arguably, this decoupling better defines a particular problem and promotes the use of diverse formulations and transformation techniques. For instance, by operating at a different level of abstraction, one might more easily identify alternative modeling and solution techniques: different measures, non-traditional support schemes, alternative derivative approximations, and/or advanced problem transformations. For instance, analyzing the problem in its infinite-dimensional form is what inspired the case study discussed in Section \ref{sec:cases_cvar}; this example shows how to use CVaR as a way to manipulate time-dependent trajectories.
\\

\begin{figure}[!htb]
    \centering
    \includegraphics[width=0.4\textwidth]{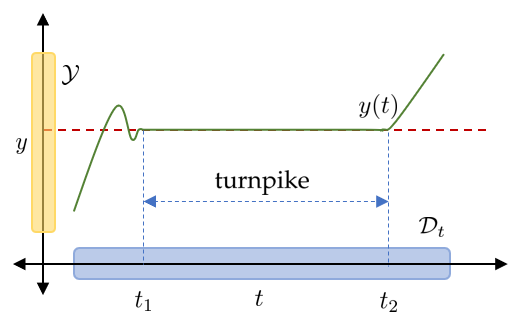}
    \caption{Illustration of the turnpike phenomenon for a time-dependent trajectory $y(t)$ where the turnpike occurs on the interval $[t_1, t_2]$.}
    \label{fig:turnpike}
\end{figure}

Establishing theoretical properties for InfiniteOpt formulations is also often facilitated when these problems are expressed in their native form. This is exemplified in the recent work of Faulwasser and Grüne \cite{faulwasser2020turnpike}, where the authors utilize continuous and discrete time formulations to derive properties of the turnpike phenomenon in the field of optimal control. The turnpike phenomenon refers to the propensity of optimal control trajectories to remain within a certain region for a significant portion of the time horizon until later departing it. Figure \ref{fig:turnpike} illustrates this for a dynamic variable $y(t)$. The properties discussed by the authors with regard to turnpike behavior are beyond the scope of this work but, interestingly, they observe that a considerable amount of analysis has been done for finite time formulations whereas many conceptual gaps remain for the continuous-time case. This conceptual disparity between the continuous and discrete time cases can at least in part be attributed to the rather standard practice of expressing optimal control formulations in discrete time. This observation is not unique to the optimal control community and there exists much to be explored for InfiniteOpt formulations in their native forms throughout their respective communities in general. Some potential avenues of research for InfiniteOpt formulations include systematic initialization techniques that consider the formulation infinite domain, generalized pre-solving methods (e.g., feasibility checking), and enhanced transformations (e.g., basis function approaches used in \cite{georgakis2013design} and \cite{gnegel2021solution}).  

\begin{figure}[!htb]
    \centering
    \begin{subfigure}[b]{0.4\textwidth}
        \centering
        \includegraphics[width=\textwidth]{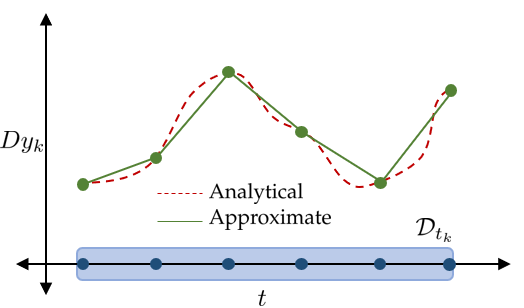}
        \caption{Traditional}
        \label{fig:bad_derivative}
    \end{subfigure}
    \quad \quad \quad
    \begin{subfigure}[b]{0.4\textwidth}
        \centering
        \includegraphics[width=\textwidth]{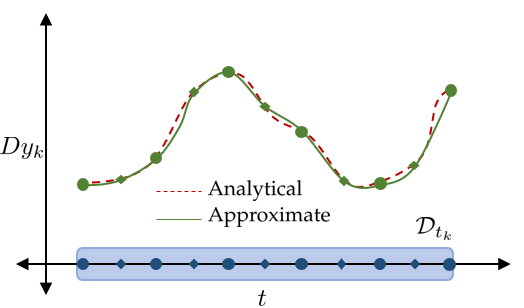}
        \caption{New}
        \label{fig:good_derivative}
    \end{subfigure}
    \caption{Comparison between the derivative approximation approaches common to traditional dynamic estimation formulations and higher-order ones possible using our new formulation (e.g., using orthogonal collocation).}
    \label{fig:derivative_estimation}
\end{figure}

The proposed abstraction also encourages a more systematic treatment of infinite-dimensional data and modeling objects. To illustrate this, we consider dynamic parameter estimation formulations and show how lifting a finite-dimensional data representation into an infinite-dimensional form might be beneficial. We consider conducting dynamic experiments $k \in \mathcal{K}$ that collect the set of observations (data) $\{\tilde{y}_k(t): t \in \hat{\mathcal{D}}_{t_k}, \ k \in \mathcal{K}\}$ over the set of time points $t \in \hat{\mathcal{D}}_{t_k}$. A dynamic parameter estimation formulation then seeks the optimal choice of parameters $z \in \mathcal{Z}$ to fit and verify the efficacy of a candidate model $g(y(t), z, t) = 0$ relative to empirical data \cite{shin2019scalable, venturelli2018deciphering}. For simplicity in example, we consider a least-squares approach that yields the following canonical discrete-time estimation formulation: 
\begin{equation}
    \begin{aligned}
       \min_{y_{k}(\cdot), z}  &&& \sum_{k \in \mathcal{K}}  \sum_{t \in \hat{\mathcal{D}}_{t_k}} (y_{k}(t)-\tilde{y}_{k}(t))^2  \\
       \text{s.t.} &&& \hat{g}(y(t), z, t) = 0, \ \ t \in \hat{\mathcal{D}}_{t_k}, k \in \mathcal{K} \\
       &&& z \in \mathcal{Z}
    \end{aligned}
    \label{eq:discrete_estimation}
\end{equation}
where $\hat{g}(\cdot)$ is the discretized dynamic model and $y_{k}(\cdot)$ are the model predicted variables. This discrete representation is guided by the nature of the experimental data $\tilde{y}_k(t)$ available, which corresponds to a finite set of time points $t \in \hat{\mathcal{D}}_{t_k}$. This limits the model to a particular transcribed domain; however, we can express formulation \eqref{eq:discrete_estimation} in continuous time by representing the experimental data with appropriate infinite-dimensional lifting functions $f_k(t), t \in \mathcal{D}_{t_k}$ such that $f_k(t)=\tilde{y}(t)$ at $t\in \hat{\mathcal{D}}_{t_k}$:
\begin{equation}
    \begin{aligned}
       \min_{y_{k}(\cdot), z}  &&& \sum_{k \in \mathcal{K}} \left( \int_{t \in \mathcal{D}_{t_k}} (y_{k}(t)-\tilde{y}_{k}(t))^2dt \right)  \\
       \text{s.t.} &&& g(y(t), z, t) = 0, &&  t \in \mathcal{D}_{t_k}, k \in \mathcal{K} \\ 
       &&& \tilde{y}_k(t) = f_k(t), &&  t \in \mathcal{D}_{t_k}, k \in \mathcal{K}\\
       &&& z \in \mathcal{Z}.
    \end{aligned}
    \label{eq:cont_estimation}
\end{equation}
We now have a formulation that fits into our unifying InfiniteOpt abstraction; as such, we can begin to consider general modeling elements and transformations. This means, for instance, that we might want to consider alternative time-dependent measures or more accurate derivative approximations (e.g., orthogonal collocation over finite elements) \cite{tjoa1991simultaneous}. Figure \ref{fig:derivative_estimation} demonstrates this principle for a certain derivative $D\,y_k$. This approach also has the potential to alleviate the large computational burden associated with the noisy experimental data that often plague dynamic estimation \cite{ramsay2007parameter}, since the chosen empirical data functions $f_k(\cdot)$ smooth the empirical domains as a preprocessing step (see Figure \ref{fig:fit_plot} in the case study). The data functions also facilitate the computation of data derivatives; which can be used in estimation techniques such as \texttt{SINDy}.  We leave a more rigorous analysis of formulation \eqref{eq:cont_estimation} to future work but we hope that this discussion helps illustrate how lifting can help identify new perspectives to tackle problems that are typically treated as finite-dimensional. We study this approach further in the biological dynamic parameter estimation case study presented in Section \ref{sec:estimate}.

\section{Case Studies} \label{sec:cases}

In this section, we provide illustrative case studies to demonstrate the concepts discussed. These case studies seek to exemplify how the unifying abstraction captures a wide range of formulation classes and how it drives innovation. These cases are implemented via \texttt{InfiniteOpt.jl v0.4.1} using \textrm{Ipopt v3.13.2} for continuous problems and \texttt{Gurobi v9.1.1} for integer-valued problems on an Intel\textregistered \, Core\texttrademark \, i7-7500U machine running at 2.90 GHz with 4 hardware threads and 16 GB of RAM running Windows 10 Home. All scripts needed to reproduce the results can be found in \url{https://github.com/zavalab/JuliaBox/tree/master/InfiniteDimensionalCases}.

\subsection{Event-Constrained Optimal Power Flow} \label{sec:sopf}

We apply the event constraints featured in Section \ref{sec:event_constrs} to a stochastic optimal power flow (SOPF) formulation (a stochastic optimization problem). We base our SOPF formulation as a variant of the chance-constrained formulation presented in \cite{muhlpfordt2018generalized}. This considers DC power grid networks subject to random power demands $\xi \in \mathcal{D}_\xi \subseteq \mathbb{R}^{n_\xi}$. The optimal policy defines the power generation $y_g(\xi) \in \mathcal{Y}_g \subseteq \mathbb{R}^{n_g}$ and the branch power flow $y_b(\xi) \in \mathcal{Y}_b \subseteq \mathbb{R}^{n_b}$ recourse functions to satisfy the demands where the respective feasible sets denote the engineering limits (i.e., $\mathcal{Y}_g = [0, \overline{y_g}]$ and $\mathcal{Y}_b = [-\overline{y_b}, \overline{y_b}]$). The power model enforces a linear energy balance at each node of the form:
\begin{equation}
    A y_b(\xi) + C_g y_g(\xi) - C_\xi \xi = 0,\; \xi \in \mathcal{D}_\xi
\end{equation}
where $A \in \mathbb{R}^{n_n \times n_b}$ is the incidence matrix, $C_g \in \mathbb{R}^{n_n \times n_g}$ maps the generators to the correct nodes, and $C_\xi \in \mathbb{R}^{n_n \times n_\xi}$ maps the demands to the correct nodes.
\\

A traditional joint-chance constraint enforces limits to a certain probabilistic threshold $\alpha$:
\begin{equation}                             
    \mathbb{P}_\xi\left(\left(\bigcap_{i \in \mathcal{I}_g}\{\forall \xi \in \mathcal{D}_\xi : y_{g, i}(\xi) \leq \overline{y_g}_i\}\right) \bigcap \left(\bigcap_{i \in \mathcal{I}_b}\{\forall \xi \in \mathcal{D}_\xi : -\overline{y_b}_i \leq y_{b, i}(\xi) \leq \overline{y_b}_i\}\right) \right)\geq \alpha
    \label{eq:sopf_chance1}
\end{equation}
where $\mathcal{I}_g$ is the set of generator indices and $\mathcal{I}_b$ is the set of branch indices. The non-negativity constraints on the generators are excluded such that they are enforced almost surely. The joint-chance constraint thus enforces that all the engineering limits are satisfied to at least a probability $\alpha$ and constraint \eqref{eq:sopf_chance1} is equivalent to:
\begin{equation}                             
\mathbb{E}_\xi \left[\mathbbm{1}\left[{\left(\bigcap_{i \in \mathcal{I}_g} y_{g, i}(\xi) \leq \overline{y_g}_i\right) \bigcap \left(\bigcap_{i \in \mathcal{I}_b} -\overline{y_b}_i \leq y_{b, i}(\xi) \leq \overline{y_b}_i\right)}\right]\right] \geq \alpha.
    \label{eq:sopf_chance2}
\end{equation}
This representation can be reformulated into algebraic constraints by introducing an infinite binary variable $y_w(\xi) \in \{0, 1\}$ and an appropriate upper-bounding constant $U \in \mathbb{R}_+$ \cite{luedtke2008sample}:
\begin{equation}
    \begin{aligned}
        &&& y_{g, i}(\xi) - \overline{y_g}_i \leq y_w(\xi) U, && \xi \in \mathcal{D}_\xi, i \in \mathcal{I}_g  \\
        &&& -y_{b, i}(\xi) - \overline{y_b}_i \leq y_w(\xi) U, &&  \xi \in \mathcal{D}_\xi, i \in \mathcal{I}_b  \\
        &&& y_{b, i}(\xi) - \overline{y_b}_i \leq y_w(\xi) U, &&  \xi \in \mathcal{D}_\xi, i \in \mathcal{I}_b \\
        &&& \mathbb{E}_\xi\left[1 - y_w(\xi)\right] \geq \alpha.
    \end{aligned}
    \label{eq:chance_eqs}
\end{equation}
Similarly, we can apply the excursion probability constraint; this enforces the probability that any engineering limit violation be no more than $1-\alpha$:
\begin{equation}                             
    \mathbb{P}_\xi\left(\left(\bigcup_{i \in \mathcal{I}_g} \{\exists \xi \in \mathcal{D}_\xi : y_{g, i}(\xi) > \overline{y_g}_i\}\right) \bigcup \left(\bigcup_{i \in \mathcal{I}_b} \{\exists \xi \in \mathcal{D}_\xi : |y_{b, i}(\xi)| > \overline{y_b}_i\}\right)\right) \leq 1-\alpha.
    \label{eq:sopf_excursion1}
\end{equation}
This constraint is equivalent to a joint-chance constraint; this becomes apparent when we reformulate constraint \eqref{eq:sopf_excursion1} as:
\begin{equation}                             
    \mathbb{E}_\xi\left[\mathbbm{1}\left[{\left(\bigcup_{i \in \mathcal{I}_g} y_{g, i}(\xi) > \overline{y_g}_i\right) \bigcup \left(\bigcup_{i \in \mathcal{I}_b} |y_{b, i}(\xi)| > \overline{y_b}_i\right)}\right]\right] \leq 1-\alpha.
    \label{eq:sopf_excursion2}
\end{equation}
and then use $y_w(\xi)$ and $U$ to obtain:
\begin{equation}
    \begin{aligned}
        &&& y_{g, i}(\xi) - \overline{y_g}_i \leq y_w(\xi) U, && \xi \in \mathcal{D}_\xi, i \in \mathcal{I}_g  \\
        &&& -y_{b, i}(\xi) - \overline{y_b}_i \leq y_w(\xi) U, && \xi \in \mathcal{D}_\xi, i \in \mathcal{I}_b  \\
        &&& y_{b, i}(\xi) - \overline{y_b}_i \leq y_w(\xi) U, && \xi \in \mathcal{D}_\xi, i \in \mathcal{I}_b \\
        &&& \mathbb{E}_\xi\left[y_w(\xi)\right] \leq 1-\alpha.
    \end{aligned}
    \label{eq:excursion_eqs2}
\end{equation}
These are equivalent to the set of constraints \eqref{eq:chance_eqs} since $\mathbb{E}_\xi\left[1 - y_w(\xi)\right] \geq \alpha$ implies $\mathbb{E}_\xi\left[y_w(\xi)\right] \leq 1-\alpha$.
\\

As an example of leveraging logical operators (e.g., $\cap$ and $\cup$) to constraint more complex regions, we consider the probability that all the generator limits being satisfied or all the branch limits being satisfied:
\begin{equation}                             
    \mathbb{E}_\xi\left[\mathbbm{1}\left[{\left(\bigcap_{i \in \mathcal{I}_g} y_{g, i}(\xi) \leq \overline{y_g}_i\right) \bigcup \left(\bigcap_{i \in \mathcal{I}_b} -\overline{y_b}_i \leq y_{b, i}(\xi) \leq \overline{y_b}_i\right)}\right]\right] \geq \alpha.
    \label{eq:sopf_alt_chance}
\end{equation}
This encapsulates a wider probabilistic region relative to that of the joint-chance constraint \eqref{eq:sopf_chance2}. This can be reformulated into a system of algebraic constraints by following the same methodology outlined for the other event constraints; however, we need to use multiple infinite binary variables $y_{w, g}(\xi), y_{w, b}(\xi), y_{w, o}(\xi)$:
\begin{equation}
    \begin{aligned}
        &&& y_{g, i}(\xi) - \overline{y_g}_i \leq y_{w, g}(\xi) U, && \xi \in \mathcal{D}_\xi, i \in \mathcal{I}_g  \\
        &&& -y_{b, i}(\xi) - \overline{y_b}_i \leq y_{w, b}(\xi) U, && \xi \in \mathcal{D}_\xi, i \in \mathcal{I}_b  \\
        &&& y_{b, i}(\xi) - \overline{y_b}_i \leq y_{w, b}(\xi) U, && \xi \in \mathcal{D}_\xi, i \in \mathcal{I}_b \\
        &&& y_{w, o}(\xi) \geq y_{w, g}(\xi) + y_{w, b}(\xi) - 1, && \xi \in \mathcal{D}_\xi \\
        &&& \mathbb{E}_\xi\left[1 - y_{w, o}(\xi)\right] \geq \alpha.
    \end{aligned}
    \label{eq:alt_eqs}
\end{equation}

\begin{figure}[!htb]
    \centering
    \includegraphics[width=0.4\textwidth]{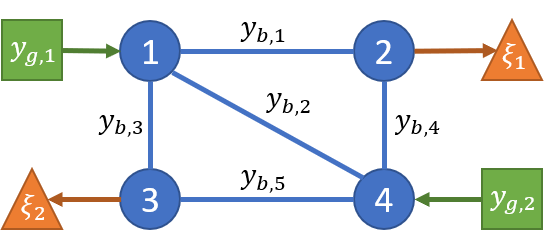}
    \caption{Sketch of the 4-bus power network topology with its bus nodes (blue circles), branches (blue lines), generators (green squares), and demand loads (orange triangles).}
    \label{fig:power_network}
\end{figure}

We can now define the SOPF formulation; here, we aim to compare the use of constraint \eqref{eq:sopf_chance2} with the use of constraint \eqref{eq:sopf_alt_chance}:
\begin{equation}
    \begin{aligned}
        &&\min_{y_g(\xi), y_b(\xi)} &&& \mathbb{E}_\xi\left[c_g^T y_g(\xi)\right] \\
        && \text{s.t.} &&& A y_b(\xi) + C_g y_g(\xi) - C_\xi \xi = 0,  && \xi \in \mathcal{D}_\xi \\
        &&&&& y_g(\xi) \geq 0, && \xi \in \mathcal{D}_\xi \\
        &&&&& \text{\eqref{eq:sopf_chance2} or \eqref{eq:sopf_alt_chance}} \\
    \end{aligned}
    \label{eq:sopf}
\end{equation}
where $c_g \in \mathbb{R}^{n_g}$ are generation unit costs. We apply this SOPF formulation to the 4-bus power system depicted in Figure \ref{fig:power_network}. We set $\overline{y_g}_i = 10$, $\overline{y_b}_i = 4$, $U = 100$, $c_g^T = [1 \, \, 10]$, and $\xi \sim \mathcal{N}(\mu, \Sigma)$ where
\begin{equation}
    \mu = \begin{bmatrix} 3 \\ 5 \end{bmatrix}, \ \ \ \Sigma = \begin{bmatrix}2 & 0 \\ 0 & 2 \end{bmatrix}.
\end{equation}
The matrices $A$, $C_g$, and $C_\xi$ are determined by the topology shown in Figure \ref{fig:power_network}. 
\\

\begin{minipage}[!htb]{0.9\linewidth}
\begin{scriptsize}
\lstset{language=Julia,breaklines = true}
\begin{lstlisting}[label = {code:infiniteopt_sopf},caption = Formulation \eqref{eq:sopf} with constraints \eqref{eq:chance_eqs} via \texttt{InfiniteOpt.jl} to obtain Pareto solutions.]
using InfiniteOpt, Distributions, Gurobi

# Initialize the model
model = InfiniteModel(Gurobi.Optimizer)

# Define the parameters
@finite_parameter(model, α == 0.95)
@infinite_parameter(model, ξ[1:2] ~ MvNormal(μ, Σ), num_supports = 1000)

# Define the variables
@variables(model, begin 
    0 <= yg[1:2], Infinite(ξ)
    yb[1:5], Infinite(ξ)
    yw, Infinite(ξ), Bin
end)

# Set the objective
@objective(model, Min, 𝔼(cg' * yg, ξ))

# Add the constraints
@constraint(model, A * yb .+ Cg * yg .- Cξ * ξ .== 0)
@constraint(model, yg - yg_lim .<= yw * U)
@constraint(model, - yb - yb_lim .<= yw * U)
@constraint(model, yb - yb_lim .<= yw * U)
@constraint(model, chance, 𝔼(1 - yw, ξ) >= α)

# Solve for the Pareto solutions 
objs = zeros(length(αs))
probs = zeros(length(αs))
for (i, prob) in enumerate(αs)
    set_value(α, prob)
    optimize!(model)
    objs[i] = objective_value(model)
    probs[i] = value(chance)
end 
\end{lstlisting}
\end{scriptsize}
\end{minipage}

\begin{figure}[!htb]
    \centering
    \includegraphics[width=0.6\textwidth]{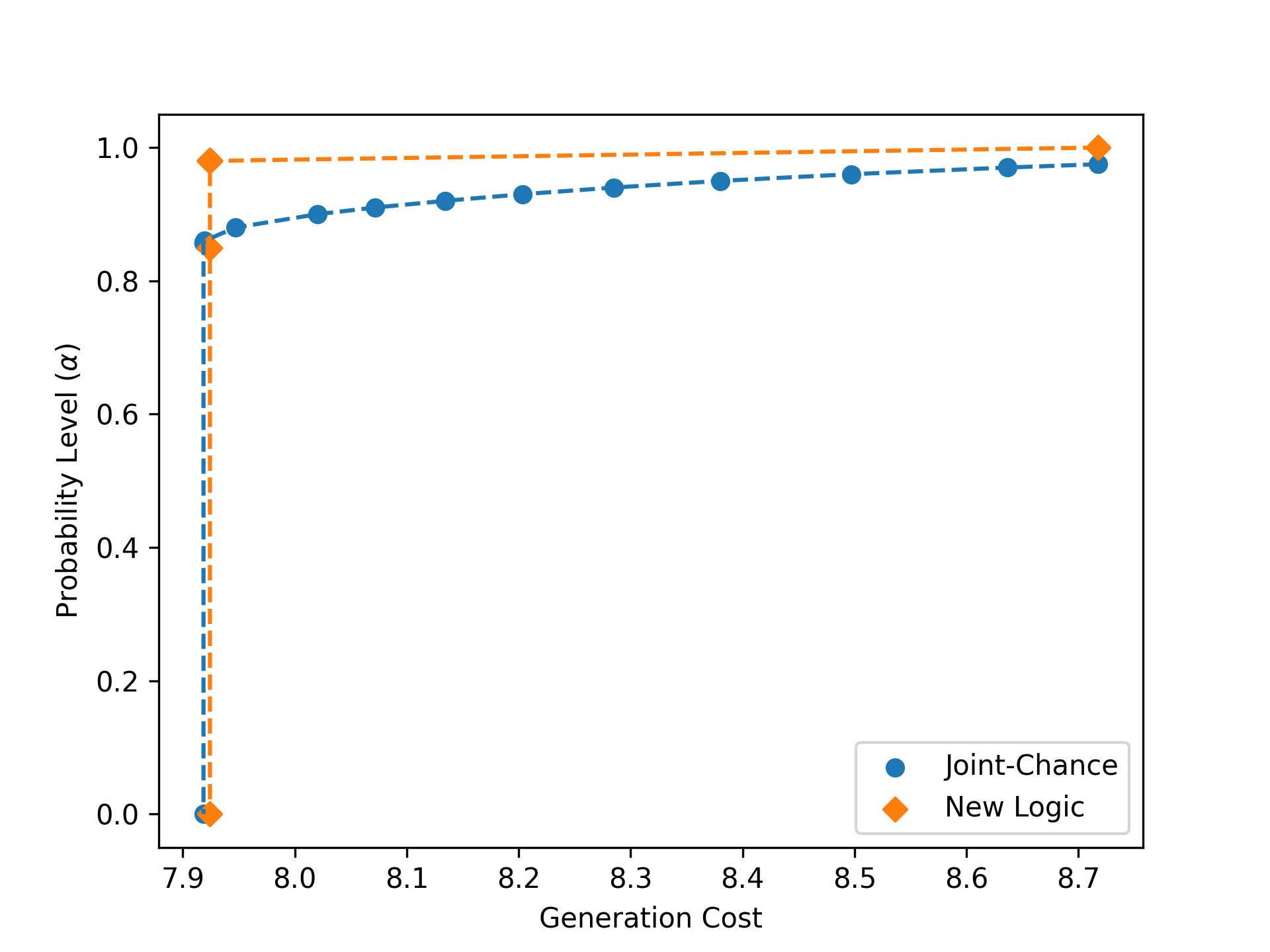}
    \caption{Pareto frontiers associated formulation \eqref{eq:sopf} where the joint-chance curve refers to using constraint \eqref{eq:sopf_chance2} and the new logic curve corresponds to constraint \eqref{eq:sopf_alt_chance}.}
    \label{fig:sopf_pareto}
\end{figure}

We implement both variants of \eqref{eq:sopf} in \texttt{InfiniteOpt.jl} and use direct transcription with 1,000 MC samples of $\xi$ to transform the InfiniteOpt problem into finite-dimensional form. Code Snippet \ref{code:infiniteopt_sopf} shows an excerpt of the script used in \texttt{InfiniteOpt.jl}. Note that the algebraic reformulations of each event constraint are used (e.g., constraints \eqref{eq:chance_eqs} in place of constraint \eqref{eq:sopf_chance2}). Each formulation is solved over a range of $\alpha$ values to obtain Pareto pairs, which are shown in Figure \ref{fig:sopf_pareto}. We observe that the traditional joint-chance constraint requires a higher power generation cost for a given probability $\alpha$. This makes sense, because the alternate formulation captures a larger probabilistic event region (it is less constraining). This highlights how logic affects the behavior of event-constrained optimization formulations. In summary, in this case study we have used our unifying abstraction to explore the use of event constraints. This is facilitated by studying the analogy between excursion probability constraints and joint-chance constraints through the lens our abstraction.

\subsection{Stochastic Optimal Pandemic Control} \label{sec:cases_cvar}

We exemplify the discussion in Section \ref{sec:risk_measures} in which we discuss how to use CVaR to shape time-dependent trajectories.  We do so by considering a pandemic optimal control problem; here, we seek to combat the spread of a contagion while minimizing the enforced isolation policy $y_u(t) \in [0, 1]$ (i.e., social distancing policy). The majority of other pandemic control studies in the literature \cite{lemecha2020optimal, area2017ebola, tsay2020modeling} use integral objectives that uniformly penalize the shape of optimal trajectories. We represent traditional formulations using the objective: 
\begin{equation}
    \min_{y_u(\cdot)} \frac{1}{S}\int_{t \in \mathcal{D}_t} y_u(t) dt
    \label{eq:pandemic_integral}
\end{equation}
with $S = \int_{t \in \mathcal{D}_t} dt$. This objective minimizes the time-average isolation policy. We also consider an alternative objective by incorporating a peak penalty (e.g., \eqref{eq:peak_cost}) to also control the maximum isolation policy:
\begin{equation}
    \min_{y_u(\cdot)} \mathop{\text{max}}_{t\in \mathcal{D}_t}\; y_u(t).
    \label{eq:pandemic_max}
\end{equation}
This new objective can be formulated as:
\begin{equation}
    \begin{aligned}
        &&\min_{y_u(\cdot), z} &&& z \\
        && \text{s.t.} &&& z \geq y_u(t),  && t \in \mathcal{D}_t \\
    \end{aligned}
    \label{eq:pandemic_max2}
\end{equation}
where $z \in \mathbb{R}$ is an auxiliary variable that captures the peak \cite{vanderbei2020linear}. We will now show that objectives \eqref{eq:pandemic_integral} and \eqref{eq:pandemic_max} are special cases of the CVaR objective:
\begin{equation}
    \min_{y_u(\cdot)} \text{CVaR}_{t}(y_u(t); \alpha)
    \label{eq:pandemic_cvar}
\end{equation}
where $\alpha \in [0, 1)$. This problem can be reformulated as:
\begin{equation}
    \begin{aligned}
        &&\min_{y_u(\cdot), y_m(\cdot), z} &&& z - \frac{1}{1-\alpha}\mathbb{E}_t[y_m(t)] \\
        && \text{s.t.} &&& y_m(t) \geq y_u(t) - z, &&  t \in \mathcal{D}_t \\
        &&&&& y_m(t) \geq 0, &&  t \in \mathcal{D}_t
    \end{aligned}
    \label{eq:pandemic_cvar2}
\end{equation}
where $y_m: \mathcal{D}_t\to \mathbb{R}$ and $z\in \mathbb{R}$ are appropriate infinite and finite auxiliary variables \cite{dowling2016framework}.

We model the spread of the contagion through a given population using the SEIR model \cite{aron1984seasonality}, which considers four population categories:
\begin{equation*}
    \text{Susceptible} \rightarrow \text{Exposed} \rightarrow \text{Infectious} \rightarrow \text{Recovered}.
\end{equation*}
We define the fractional populations of individuals susceptible to infection $y_s: \mathcal{D}_t\to [0, 1]$, exposed individuals that are not yet infectious $y_e: \mathcal{D}_t \to [0, 1]$, infectious individuals $y_i: \mathcal{D}_t\to [0, 1]$, and recovered individuals $y_r: \mathcal{D}_t\to [0, 1]$ (considered immune to future infection). The variables are normalized such that $y_s(t) + y_e(t) + y_i(t) + y_r(t) = 1$. The deterministic SEIR model is formalized as:
\begin{equation}
    \begin{gathered}
        \frac{dy_s(t)}{dt} = (y_u(t) - 1)\beta y_s(t)y_i(t),\; t\in \mathcal{D}_t \\
        \frac{dy_e(t)}{dt} = (1 - y_u(t))\beta y_s(t)y_i(t) - \xi y_e(t),\; t\in \mathcal{D}_t \\
        \frac{dy_i(t)}{dt} = \xi y_e(t) - \gamma y_i(t),\; t\in \mathcal{D}_t\\
        \frac{dy_r(t)}{dt} = \gamma y_i(t),\; t\in \mathcal{D}_t,
    \end{gathered}
\end{equation}
where $\beta, \gamma, \xi \in \mathbb{R}$ are the rates of infection, recovery, and incubation, respectively. For our case study, we consider $\xi$ to be an uncertain parameter $\xi \sim \mathcal{U}(\underline{\xi}, \overline{\xi})$. This introduces the random domain $\mathcal{D}_\xi$ (i.e., the co-domain of $\mathcal{U}(\underline{\xi}, \overline{\xi})$) and gives a stochastic dynamic optimization problem of the form:  
\begin{equation}
    \begin{aligned}
        &&\min_{} &&& \text{Objective \eqref{eq:pandemic_integral}, \eqref{eq:pandemic_max2}, or \eqref{eq:pandemic_cvar2}} \\
        && \text{s.t.} &&& \frac{\partial y_s(t, \xi)}{\partial t} = (y_u(t) - 1)\beta y_s(t, \xi) y_i(t, \xi), &&  t \in \mathcal{D}_{t}, \xi \in \mathcal{D}_{\xi} \\
        &&&&& \frac{\partial y_e(t, \xi)}{\partial t} = (1 - y_u(t))\beta y_s(t, \xi) y_i(t, \xi) - \xi y_e(t, \xi), &&  t \in \mathcal{D}_{t}, \xi \in \mathcal{D}_{\xi} \\
        &&&&& \frac{\partial y_i(t, \xi)}{\partial t} = \xi y_e(t, \xi) - \gamma y_i(t, \xi), && t \in \mathcal{D}_{t}, \xi \in \mathcal{D}_{\xi} \\
        &&&&& \frac{\partial y_r(t, \xi)}{\partial t} = \gamma y_i(t, \xi), && t \in \mathcal{D}_{t}, \xi \in \mathcal{D}_{\xi} \\
        &&&&& y_s(0, \xi) = s_0, y_e(0, \xi) = e_0, y_i(0, \xi) = i_0, y_r(0, \xi) = r_0, &&  \xi \in \mathcal{D}_{\xi} \\
        &&&&& y_i(t, \xi) \leq i_{max}, && t \in \mathcal{D}_{t}, \xi \in \mathcal{D}_{\xi} \\
        &&&&& y_u(t) \in [0, \overline{y_u}], && t \in \mathcal{D}_{t} \\
    \end{aligned}
    \label{eq:pandemic_problem}
\end{equation}
where $s_0, e_0, i_0, r_0 \in \mathbb{R}$ denote the initial population fractions, $i_{max}$ denotes the maximum allowable fraction of infected individuals $y_i(t)$, and $\overline{y_u}$ denotes the maximum realizable population isolation. The state variables $y_s(\cdot), y_i(\cdot), y_e(\cdot), y_r(\cdot)$ are now infinite variables that are parameterized in the time and random domains, while the control variable $y_u$ is an infinite variable that is only parameterized in the time domain (since we need to decide our control policy before knowing the realizations of $\xi$).

\begin{table}[!htb]
    \centering
    \caption{Parameter values used in the InfiniteOpt formulation \eqref{eq:pandemic_problem}.}
    \begin{tabular}{|c|c|c|c|c|c|c|c|c|c|}
    \hline
     $\beta$ & $\gamma$ & $\underline{\xi}$ & $\overline{\xi}$ & $i_{max}$ & $\overline{y_u}$ & $s_0$          & $e_0$     & $i_0$ & $r_0$\\ \hline \hline
     0.727   & 0.303    & 0.1               & 0.6              & 0.02      & 0.8              &  $1 - 10^{-5}$ & $10^{-5}$ & 0     & 0\\
     \hline
    \end{tabular}
    \label{tab:pandemic_params}
\end{table}

\begin{figure}[!htb]
     \centering
     \begin{subfigure}[b]{0.45\textwidth}
         \centering
         \includegraphics[width=\textwidth]{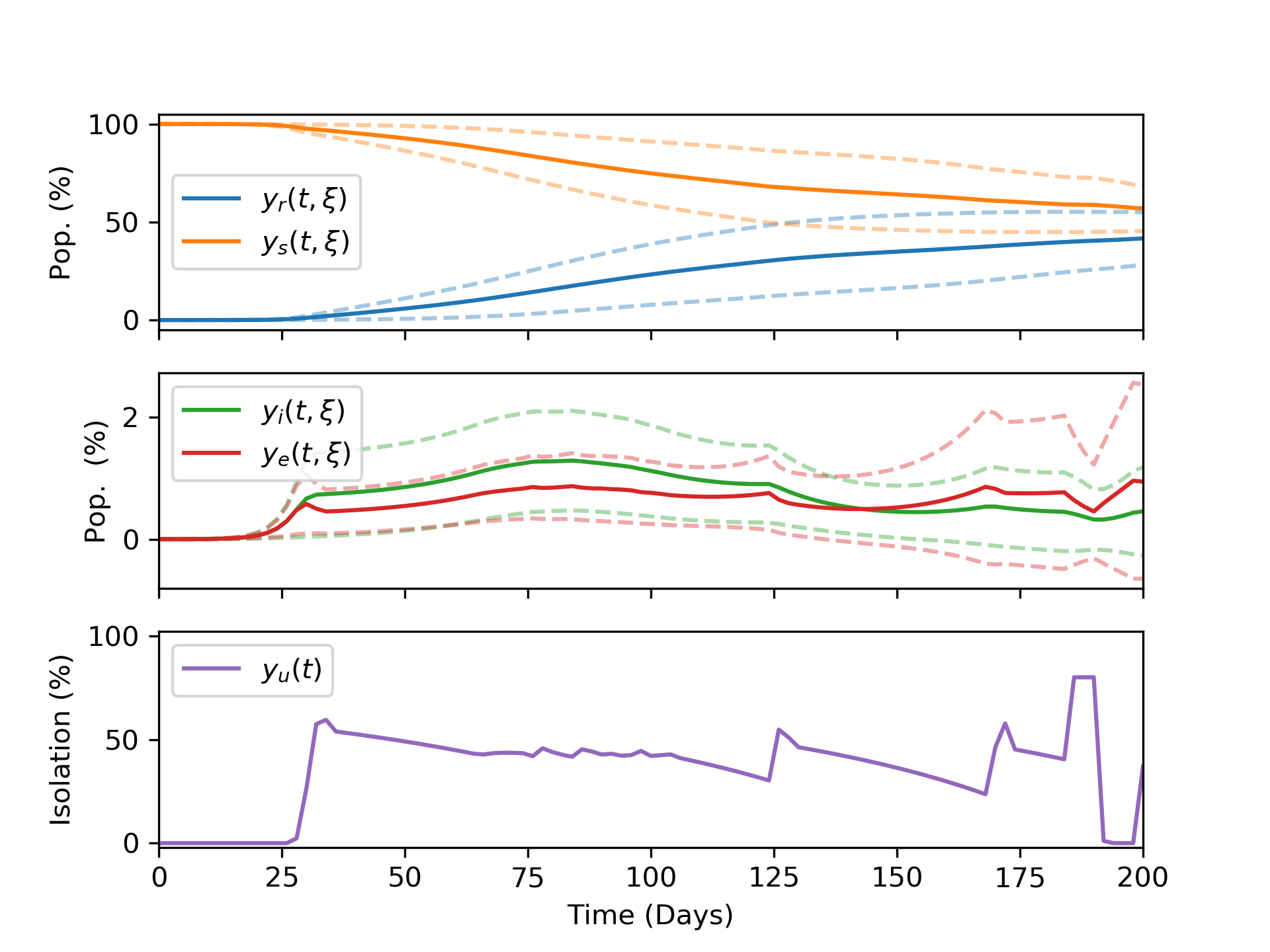}
         \caption{Objective \eqref{eq:pandemic_integral}}
         \label{fig:covid_integral}
     \end{subfigure}
     \hfill
     \begin{subfigure}[b]{0.45\textwidth}
         \centering
         \includegraphics[width=\textwidth]{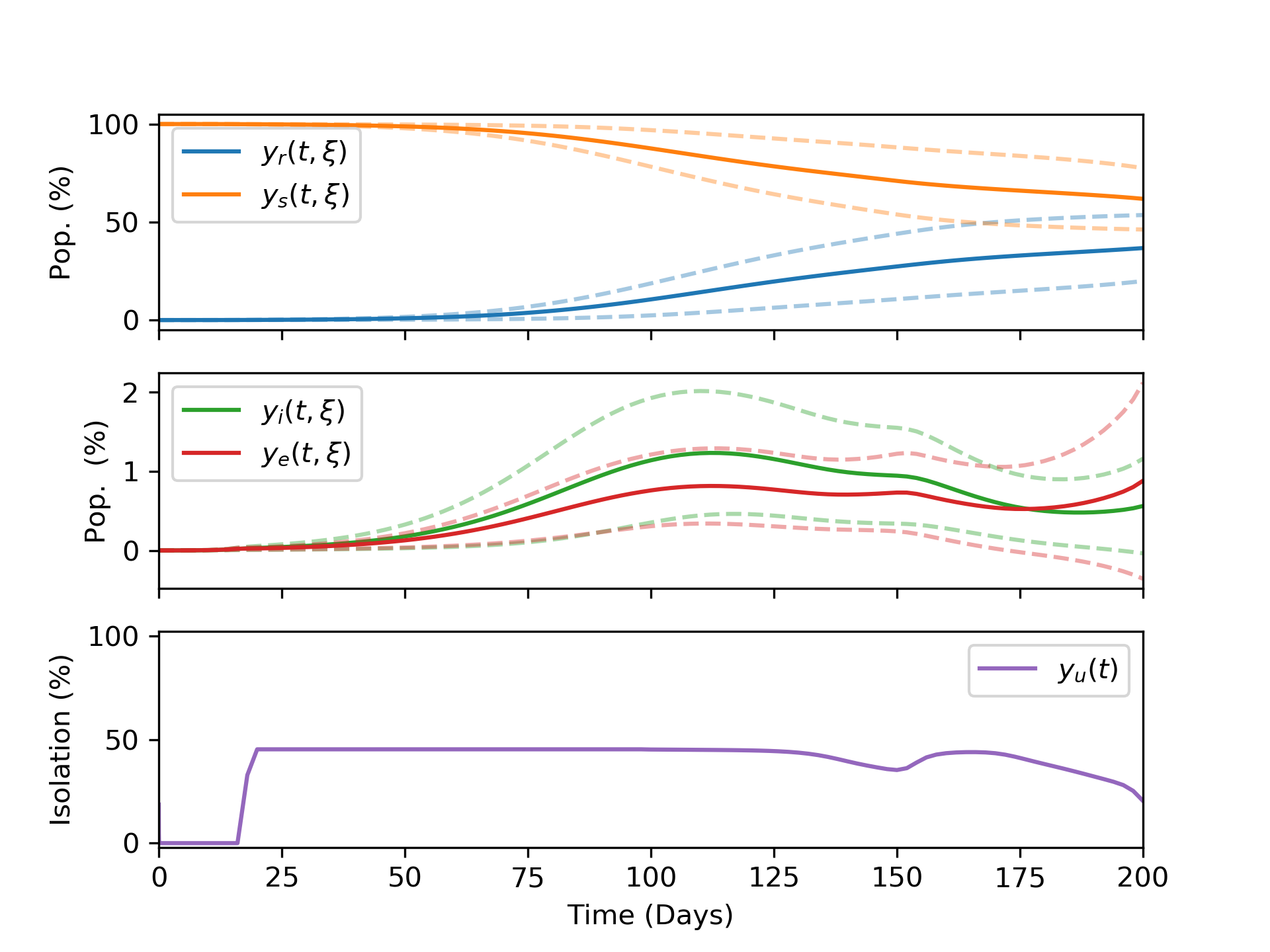}
         \caption{Objective \eqref{eq:pandemic_max2}}
         \label{fig:covid_max}
     \end{subfigure}
        \caption{Optimal trajectories for formulation \eqref{eq:pandemic_problem} using traditional dynamic optimization objectives. For the state variables $y_s(t, \xi)$, $y_e(t, \xi)$, $y_r(t, \xi)$, and $y_r(t, \xi)$ the solid lines denote the trajectories averaged over $\xi$ and the dashed lines denote the trajectories that are one standard deviation away from the mean.}
        \label{fig:standard_covid}
\end{figure}

We solve the InfiniteOpt problem \eqref{eq:pandemic_problem} using the parameters defined in Table \ref{tab:pandemic_params} with $\mathcal{D}_t = [0, 200]$. We transcribe it via \texttt{InfiniteOpt.jl} using 111 supports for $\mathcal{D}_t$ and 20 MC samples for $\mathcal{D}_\xi$. Code Snippet \ref{code:infiniteopt_covid} shows an excerpt of this implementation in \texttt{InfiniteOpt.jl}. The optimal policies corresponding to objectives \eqref{eq:pandemic_integral} and \eqref{eq:pandemic_max2} are shown in Figure \ref{fig:standard_covid}. Penalizing the peak isolation provides a smoother isolation policy $y_u(t)$ relative to the more traditional integral based objective. Moreover, the population of susceptible individuals $y_s(t, \xi)$ associated with penalizing the peak isolation decreases at a lower rate which indicates that penalizing the peak isolation is more effective at mitigating the spread of the contagion in this particular case. 

\begin{minipage}[!htb]{0.9\linewidth}
\begin{scriptsize}
\lstset{language=Julia,breaklines = true}
\begin{lstlisting}[label = {code:infiniteopt_covid},caption = Formulation \eqref{eq:pandemic_problem} with objective \eqref{eq:pandemic_cvar2} via \texttt{InfiniteOpt.jl}.]
using InfiniteOpt, Distributions, Ipopt

# Initialize the model
model = InfiniteModel(Ipopt.Optimizer)

# Set the infinite parameters 
@infinite_parameter(model, t ∈ [t0, tf], num_supports = 101)
add_supports(t, extra_ts)
@infinite_parameter(model, ξ ~ Uniform(ξ_min, ξ_max), num_supports = 20)

# Set the infinite variables 
var_inds = [:s, :e, :i, :r]
@variable(model, 0 ≤ y[var_inds], Infinite(t, ξ))
@variable(model, ysi, Infinite(t, ξ))
@variable(model, 0 ≤ yu ≤ 0.8, Infinite(t), start = 0.2)

# Set the CVaR objective 
@variable(model, z)
@variable(model, 0 ≤ ym, Infinite(t))
@objective(model, Min, z + 1 / (1 - α) * 𝔼(ym, t))
@constraint(model, ym ≥ yu - z)

# Define the initial conditions
@constraint(model, [v ∈ var_inds], y[v](0, ξ) == y0[v])

# Define the SEIR equations
@constraints(model, begin 
    ∂(y[:s], t) == -(1 - yu) * β * ysi
    ∂(y[:e], t) == (1 - yu) * β * ysi - ξ * y[:e]
    ∂(y[:i], t) == ξ * y[:e] - γ * y[:i]
    ∂(y[:r], t) == γ * y[:i]
    ysi == y[:s] * y[:i]
end)

# Define the infection limit
@constraint(model, y[:i] ≤ i_max)

# Optimize and get the results
optimize!(model)
state_opt = value.(y, ndarray = true)
control_opt = value(yu)
obj_opt = objective_value(model)
ts = value(t)
ξs = value(ξ)
\end{lstlisting}
\end{scriptsize}
\end{minipage}

\begin{figure}[!htb]
     \centering
     \begin{subfigure}[b]{0.45\textwidth}
         \centering
         \includegraphics[width=\textwidth]{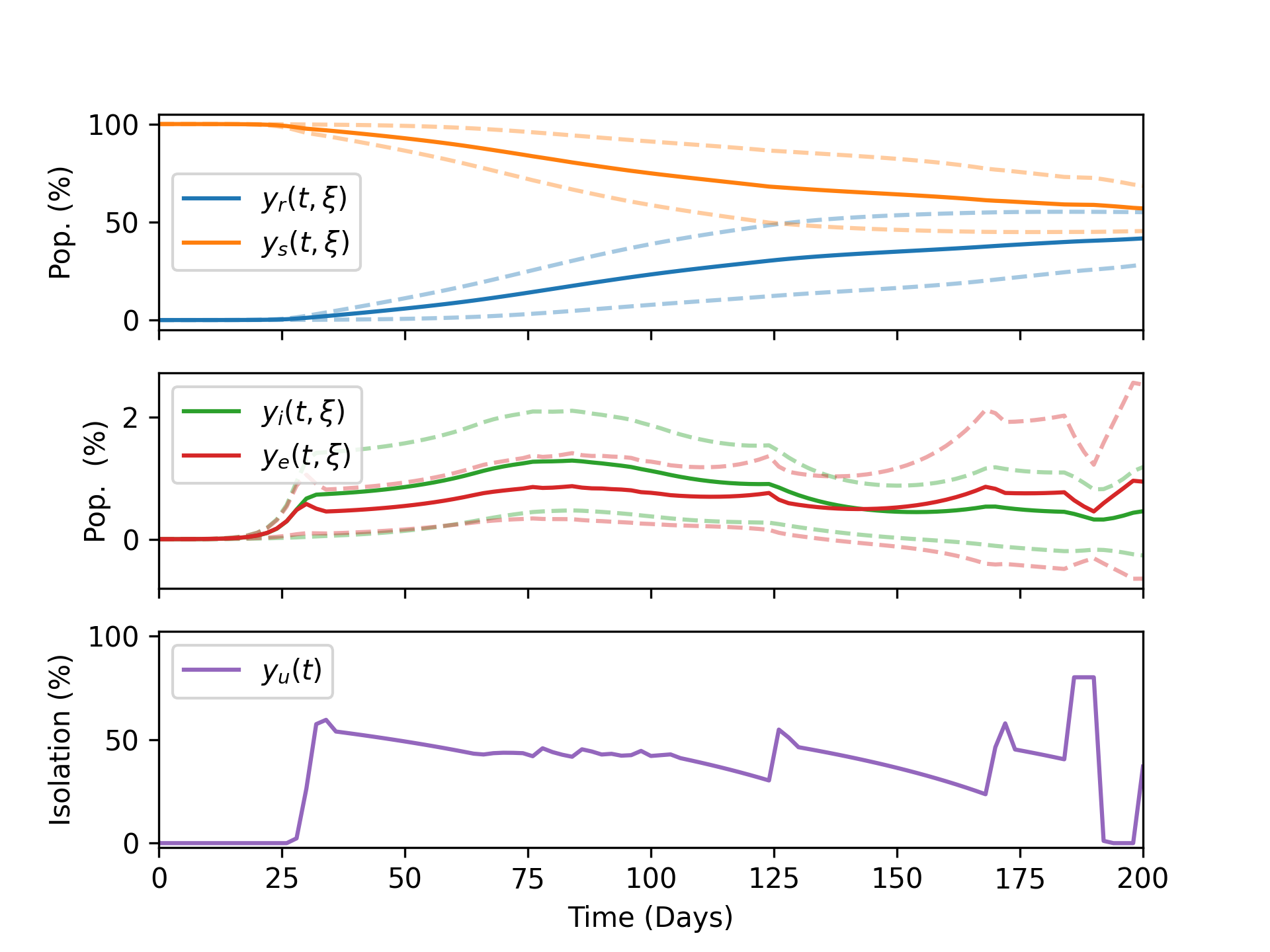}
         \caption{$\alpha = 0$}
         \label{fig:cvar_0}
     \end{subfigure}
     \hfill
     \begin{subfigure}[b]{0.45\textwidth}
         \centering
         \includegraphics[width=\textwidth]{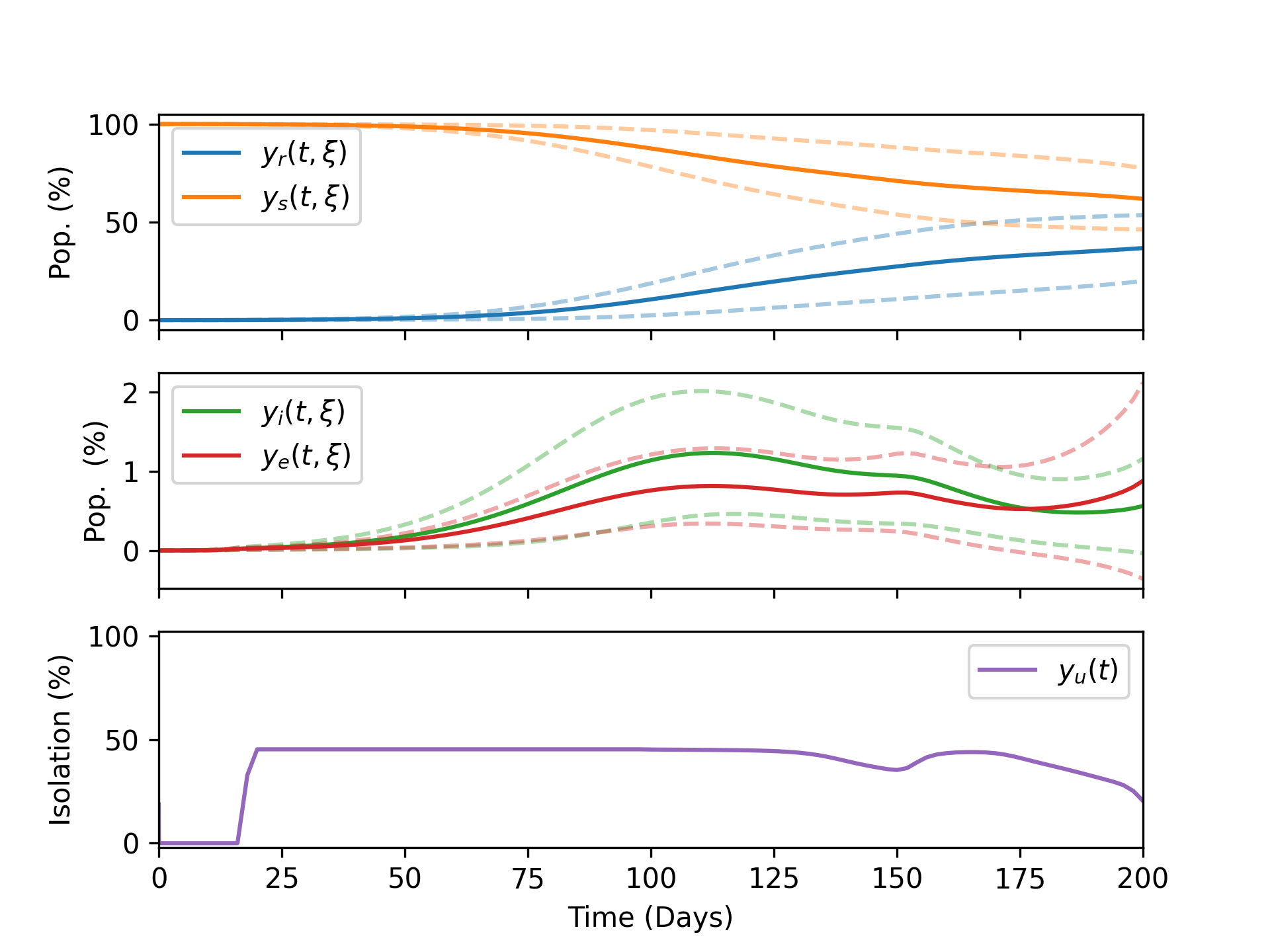}
         \caption{$\alpha = 0.99$}
         \label{fig:cvar_1}
     \end{subfigure} \\
     \begin{subfigure}[b]{0.45\textwidth}
         \centering
         \includegraphics[width=\textwidth]{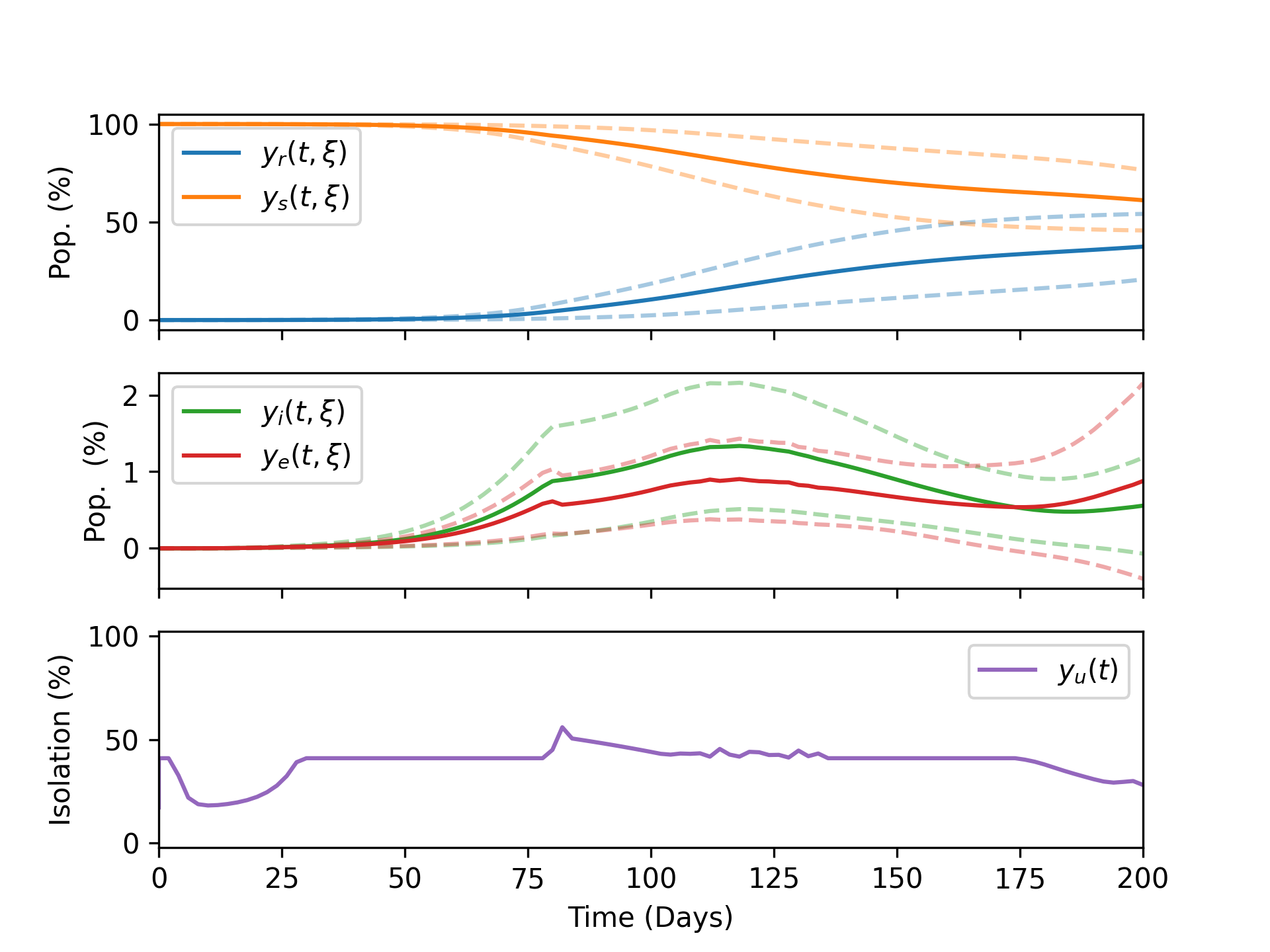}
         \caption{$\alpha = 0.5$}
         \label{fig:cvar_0_5}
     \end{subfigure}
        \caption{Optimal policies for solving Formulation \eqref{eq:pandemic_problem} in combination with Objective \eqref{eq:pandemic_cvar2}.}
        \label{fig:cvar_covid}
\end{figure}

We address formulation \eqref{eq:pandemic_problem} by using the proposed CVaR measure. Three solutions are obtained corresponding to $\alpha = \{0, 0.5, 0.99\}$ and are depicted in Figure \ref{fig:cvar_covid}. Note that we use $\alpha = 0.99$ because a singularity is incurred with $\alpha = 1$. We observe that Figures \ref{fig:cvar_0} and \ref{fig:cvar_1} are identical to Figures \ref{fig:covid_integral} and \ref{fig:covid_max}; this illustrates that the CVaR objective \eqref{eq:pandemic_cvar} is a generalization of objectives \eqref{eq:pandemic_integral} and \eqref{eq:pandemic_max}. Moreover, CVaR provides a range between these cases by varying the value of $\alpha$. Hence, we can also seek tradeoff solutions such as what is shown in Figure \eqref{fig:cvar_0_5} with $\alpha = 0.5$. Interestingly, this optimal policy combines aspects of the average and peak solution policies, but also is unique in enforcing isolation policies from the start.  This shows that CVaR can be used to shape time-dependent trajectories in unique ways that are difficult to achieve with traditional measures used in dynamic optimization.

\subsection{Estimation in Dynamic Biological Systems} \label{sec:estimate}

We consider a dynamic parameter estimation problem for a biological system model; this study aims to demonstrate how the unifying abstraction inspires new formulation classes by lifting formulations into infinite dimensional spaces. This follows from the discussion in Section \ref{sec:characterizations} with regard to formulations \eqref{eq:discrete_estimation} and \eqref{eq:cont_estimation}.

\begin{table}[!htb]
    \centering
    \caption{Species membership of the microbial community.}
    \begin{tabular}{|c | c|} 
        \hline
        Species Name & Abbreviation \\ [0.5ex] 
        \hline \hline
        Blautia hydrogenotrophica & BH \\ \hline
        Collinsella aerofaciens & CA \\ \hline
        Bacteroides uniformis & BU \\ \hline
        Prevotella copri & PC \\ \hline
        Bacteroides ovatus & BO \\ \hline
        Bacteroides vulgatus & BV \\ \hline
        Bacteroides thetaiotaomicron & BT \\ \hline
        Eggerthella lenta & EL \\ \hline
        Faecalibacterium prausnitzii & FP \\ \hline
        Clostridium hiranonis & CH \\  \hline
        Desulfovibrio piger & DP \\ \hline
        Eubacterium rectale & ER \\
        \hline
    \end{tabular}
    \label{tab:microbe_names}
\end{table}

To juxtapose the utility of formulations \eqref{eq:discrete_estimation} and \eqref{eq:cont_estimation}, we consider a microbial community consisting of the 12 species described in Table \ref{tab:microbe_names} using the generalized Lotka-Volterra (gLV) model:
\begin{equation}
     \frac{dy_{x,i}(t)}{dt} = \bigg(z_{\mu,i} + \sum_{j \in \mathcal{I}} z_{\alpha, ij}y_{x,j}(t)\bigg) y_{x,i}(t),\; \ i \in \mathcal{I}
     \label{eq:glv}
\end{equation}
where $\mathcal{I}$ represents the set of microbial species, $y_{x,i}(t)$ is the estimated absolute abundance of species $i \in \mathcal{I}$, $z_{\mu,i}$ is the growth rate of species $i$, and $z_{\alpha, ij}$ is a species interaction coefficient which describes the effect of the recipient species $i$ on the growth of the donor species $j$ \cite{shin2019scalable}. We use the gLV model parameters presented in \cite{venturelli2018deciphering} to generate simulated experimental data with random noise $\epsilon \sim \mathcal{N}(0, 0.01)$ for 12 mono-species and 66 pairwise experiments. This will enhance our assessment of Formulations \eqref{eq:discrete_estimation} and \eqref{eq:cont_estimation}, since we have an established ground truth for the model parameters. 
\\

Incorporating \eqref{eq:glv} into Formulation \eqref{eq:discrete_estimation} provides us with our discrete dynamic biological community estimation formulation:
\begin{equation}
    \begin{aligned}
     \min  &&& \sum_{k \in \mathcal{K}} \sum_{i \in \mathcal{I}} \sum_{t\in \hat{\mathcal{D}}_{t_k}} (y_{x,ik}(t)-\tilde{y}_{x, ik}(t))^2 \\
    \text{s.t.} &&& \frac{dy_{x,ik}(t)}{dt} = \bigg(z_{\mu,i} + \sum_{j \in \mathcal{I}} z_{\alpha, ij}y_{x,jk}(t)\bigg) y_{x,ik}(t), &&  t \in \hat{\mathcal{D}}_{t_k}, i \in \mathcal{I}, k \in \mathcal{K} \\
    &&& 0.09 \le z_{\mu, i} \le 2.1, &&  i \in \mathcal{I} \\
    &&& -10 \le z_{\alpha, ii} \le 0, &&  i \in \mathcal{I} \\
    &&& -10 \le z_{\alpha, ij} \le 10, && (i, j \neq i) \in \mathcal{I} \times \mathcal{I}.
    \end{aligned}
    \label{eq:full_discrete_estimate}
\end{equation}
Note that we sum over each species $i$ for each experiment $k$ in accordance with \eqref{eq:glv}. Also, we recall that the derivative terms are limited to finite difference approximation schemes that employ the support points in $\hat{\mathcal{D}}_{t_k}$.

\begin{figure}[!htb]
    \centering
    \includegraphics[width=0.5\textwidth]{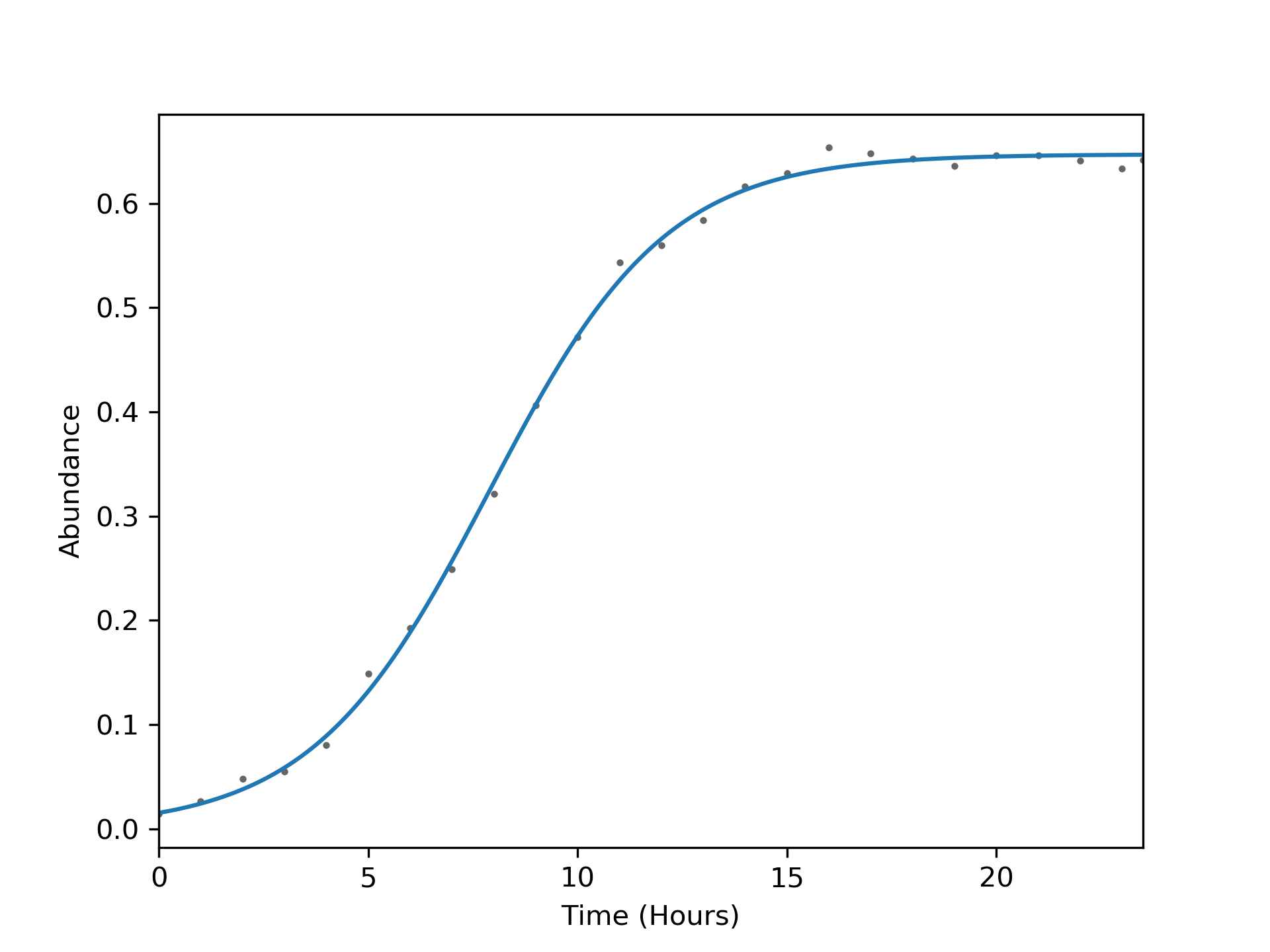}
    \caption{Empirical fit for a mono-species experiment of Bacteroides vulgatus using \eqref{eq:sigmoid_func}.}
    \label{fig:fit_plot}
\end{figure}

To derive an explicit form of Formulation \eqref{eq:cont_estimation}, we first fit an empirical function $f_{ik}(t_k)$ to each experiment $k$ and each species $i$. Our goal is not to create a generalizable predictive representation but rather to construct a continuous, infinite-dimensional function that represents the dynamic behavior observed in each experiment. We observe that the experiments appear to exhibit sigmoidal characteristics, and thus we use least-squares to fit each experiment $k$ and species $i$ to an empirical function of the form: 
\begin{equation}
    f_{ik}(t) := \frac{\beta_{1,ik}}{\beta_{2,ik} + \beta_{3,ik} e^{\beta_{4,ik}(t - \beta_{5,ik})}},\; t \in \mathcal{D}_{t_k}, i \in \mathcal{I}, k \in \mathcal{K}
    \label{eq:sigmoid_func}
\end{equation}
where $\beta_{1,ik}$, $\beta_{2,ik}$, $\beta_{3,ik}$, $\beta_{4,ik}$, and $\beta_{5,ik}$ are fitting parameters. Equation \eqref{eq:sigmoid_func} fits our experimental datasets well as demonstrated for the particular experiment shown in Figure \ref{fig:fit_plot} which is indicative of overall fit qualities we observed. Thus, we substitute Equations \eqref{eq:glv} and \eqref{eq:sigmoid_func} into Formulation \eqref{eq:cont_estimation} to obtain:
\begin{equation}
    \begin{aligned}
     \min  &&& \sum_{k \in \mathcal{K}} \sum_{i \in \mathcal{I}} \bigg(\int_{t \in \mathcal{D}_{t_k}} (y_{x,ik}(t)-\tilde{y}_{x, ik}(t))^2dt\bigg) \\
     \text{s.t.} &&& \frac{dy_{x,ik}(t)}{dt} = \bigg(z_{\mu,i} + \sum_{j \in \mathcal{I}} z_{\alpha, ij}y_{x,jk}(t)\bigg) y_{x,ik}(t), && t \in \mathcal{D}_{t_k}, i \in \mathcal{I}, k \in \mathcal{K} \\
    &&& \tilde{y}_{x, ik}(t) = \frac{\beta_{1,ik}}{\beta_{2,ik} + \beta_{3,ik} e^{\beta_{4,ik}(t - \beta_{5,ik})}}, && t \in \mathcal{D}_{t_k}, i \in \mathcal{I}, k \in \mathcal{K} \\
    &&& 0.09 \le z_{\mu, i} \le 2.1, &&  i \in \mathcal{I} \\
    &&& -10 \le z_{\alpha, ii} \le 0, &&  i \in \mathcal{I} \\
    &&& -10 \le z_{\alpha, ij} \le 10, && (i, j \neq i) \in \mathcal{I} \times \mathcal{I}.
    \end{aligned}
    \label{eq:full_cont_estimate}
\end{equation}

\begin{minipage}[!htb]{0.9\linewidth}
\begin{scriptsize}
\lstset{language=Julia,breaklines = true}
\begin{lstlisting}[label = {code:infiniteopt_estimate},caption = Formulation \eqref{eq:full_cont_estimate} using \texttt{InfiniteOpt.jl}.]
using InfiniteOpt, Ipopt

# Initialize the model
model = InfiniteModel(Ipopt.Optimizer) 

# Set the infinite parameters
@infinite_parameter(model, t[k ∈ K] ∈ [0, T[k]], num_supports = 15, independent = true, 
                    derivative_method = OrthogonalCollocation(4))

# Set the finite variables
@variable(model, 0.09 ≤ zμ[i ∈ I] ≤ 2.1)
@variable(model, -10 ≤ zα[i ∈ I, j ∈ I] ≤ zα_max[i, j])

# Set the infinite variables
@variable(model, yx[i ∈ I, k ∈ K] ≥ 0, Infinite(t[k]))
@variable(model, zαyx[i ∈ I, j ∈ I, k ∈ K], Infinite(t[k]))

# Set the empirical functions using fitted q[i, k](t) functions
@parameter_function(model, yx_tilde[i ∈ I, k ∈ K] == q[i, k](t[k]))

# Set the least-squares objective 
@objective(model, Min, sum(∫((yx[i, k] - yx_tilde[i, k]) ^ 2, t[k]) for i ∈ I, k ∈ K))

# Define the gLV equations
@constraint(model, [i ∈ I, k ∈ K], 
            ∂(yx[i, k], t[k]) == (zμ[i] + sum(zαyx[i, j, k] for j ∈ I)) * yx[i, k])
@constraint(model, [i ∈ I, j ∈ I, k ∈ K], zαyx[i, j, k] == zα[i, j] * yx[j, k])

# Optimize and get the results
optimize!(model)
zα_opt = value.(zα)
zμ_opt = value.(zμ)
yx_opt = value.(yx)
ts = supports(t)
\end{lstlisting}
\end{scriptsize}
\end{minipage}
    
\begin{figure}[!htb]
    \centering
    \includegraphics[width=1.0\textwidth]{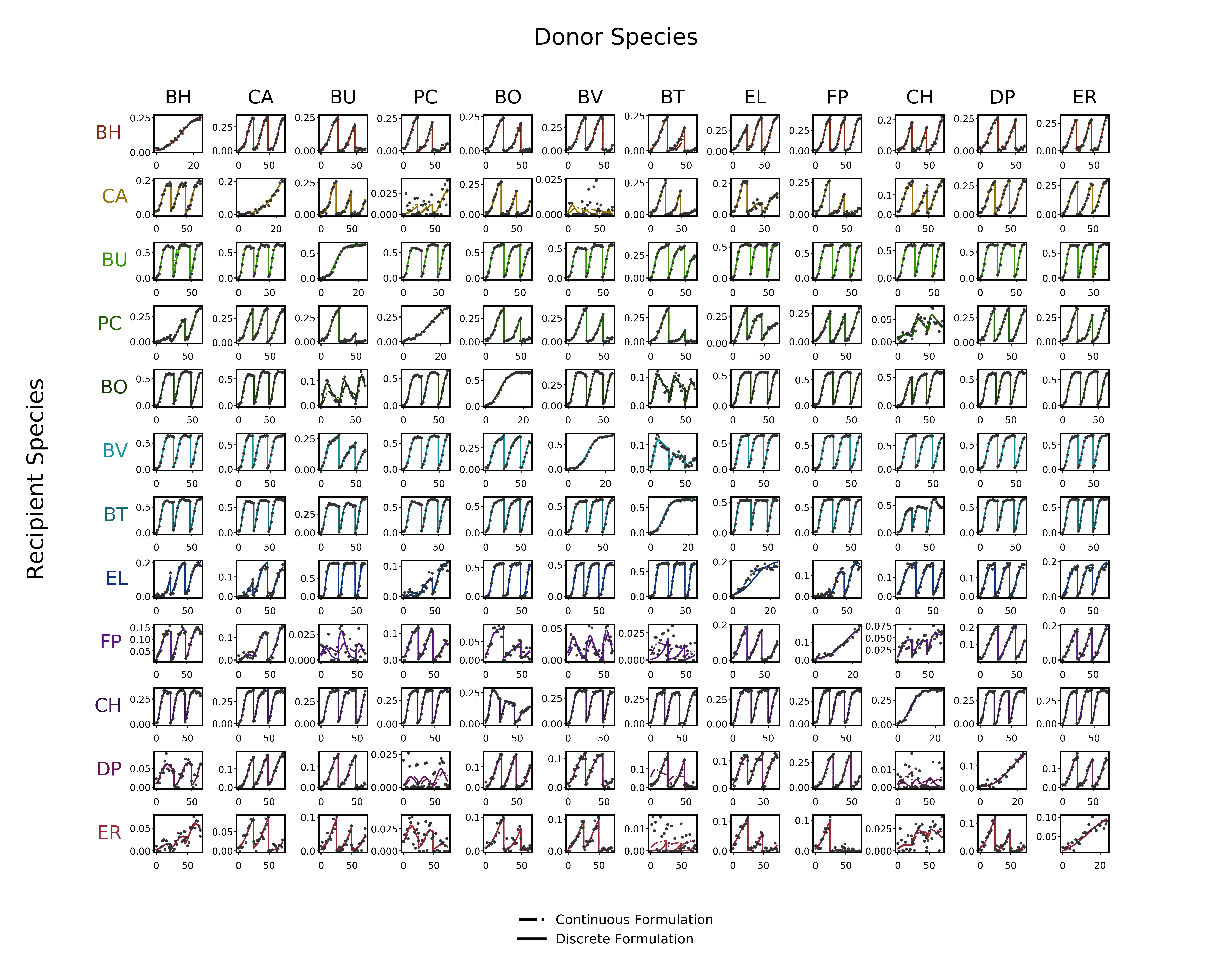}
    \caption{Optimal trajectories from formulations \eqref{eq:full_discrete_estimate} and \eqref{eq:full_cont_estimate} using orthogonal collocation over finite elements to approximate the derivatives with two and four points, respectively. All shown in comparison to the experimental data. The x-axis is time in hours, and the y-axis is the absolute abundance of the recipient species in contact with the corresponding donor species. The results for the mono-species experiments are observed on the diagonal with the rest being pairwise.}
    \label{fig:all_fits}
\end{figure}

We solve formulations \eqref{eq:full_discrete_estimate} and \eqref{eq:full_cont_estimate} using \texttt{InfiniteOpt.jl} via its automated transcription capabilities. Code Snippet \ref{code:infiniteopt_estimate} presents an illustrative summary of the syntax for implementing Formulation \eqref{eq:full_cont_estimate}. We employ 15 time supports for each experiment in Formulation \eqref{eq:full_cont_estimate} and necessarily use the measurement times as supports for formulation \eqref{eq:full_discrete_estimate}. We use orthogonal collocation over finite elements to approximate the differential equations; we use two nodes per finite element (necessarily using only the boundary points) in Formulation \eqref{eq:full_discrete_estimate} and a range of node amounts with formulation \eqref{eq:full_cont_estimate} to investigate their effect on solution quality. Figure \ref{fig:all_fits} summarizes the model fits of both solution sets relative to the experimental data. Although both estimation problems are able to choose parameters that characterize the data well, there are significant deviations between the profiles across a few experiments. 

\begin{figure}[!htb]
    \centering
    \includegraphics[width=0.8\textwidth]{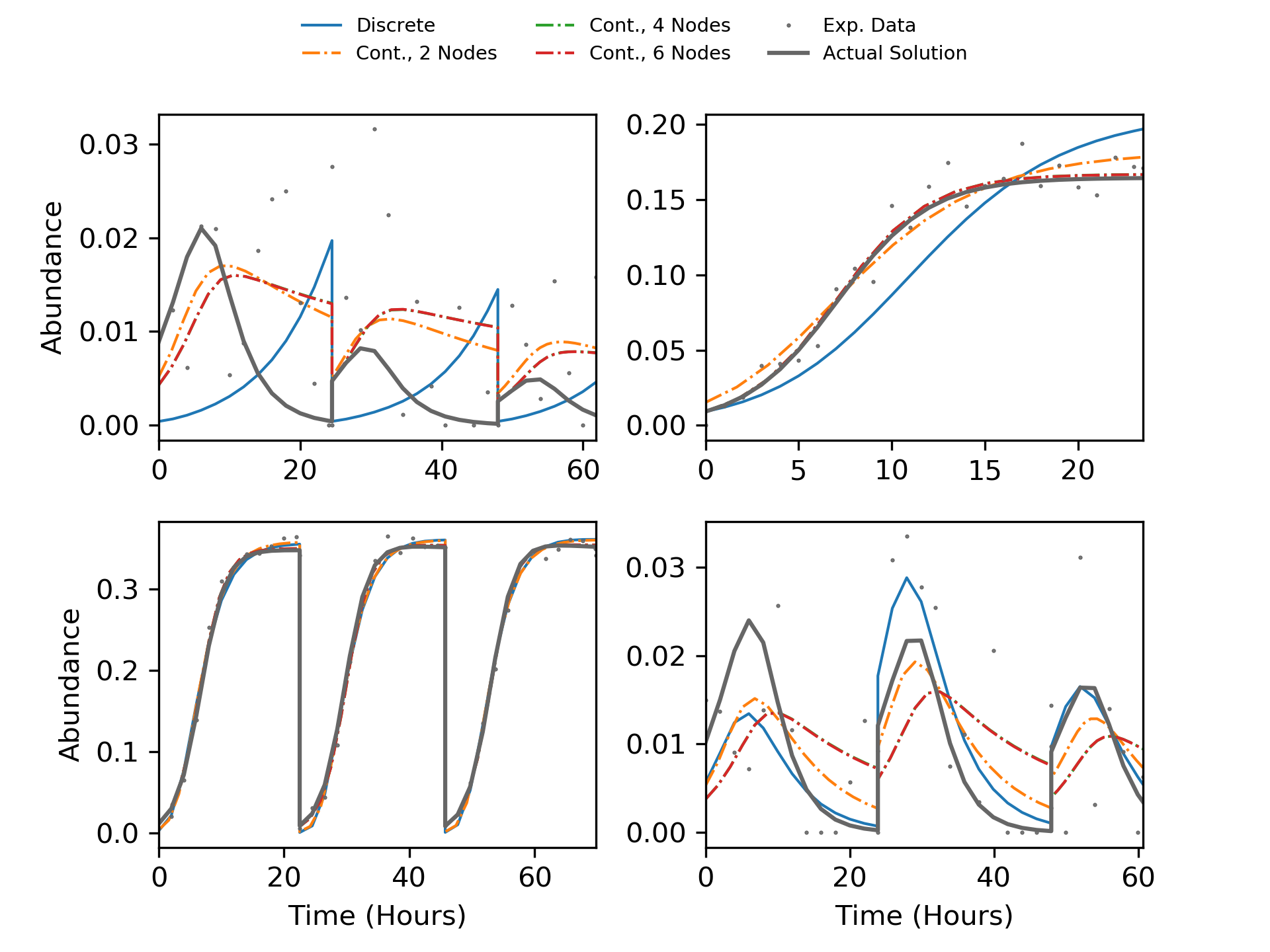}
    \caption{Optimal profiles for select experiments using different formulations and collocation node amounts (top-left: (FP, BT), top-right: (EL, EL), bottom-left: (CH, PC), bottom-right: (FP, BU)).}
    \label{fig:fit_accuracies}
\end{figure}

Figure \ref{fig:fit_accuracies} shows a representative subset of experiments to demonstrate how the solutions vary with respect to the estimation formulation and the number of collocation points. We find that the discrete formulation solution generally deviates from the analytical solution to a greater extent than the solutions procured via the continuous formulation. However, it is difficult to conclude which estimation problem best represents the data when the measurement noise is near the magnitude of the absolute abundance because no model matches the true analytical solution particularly well; the continuous formulation solutions, however, in general seem to better represent the trend of the system. This suggests that Formulation \eqref{eq:full_cont_estimate} has effectively smoothed over noisy experiments, leading to a better fit. Furthermore, we observe that the ability of Formulation \eqref{eq:full_cont_estimate} to enable arbitrary collocation nodes has a significant effect on the accuracy of the solutions. This increased accuracy seems to effectively taper off at four collocation nodes in this case. 

We seek to substantiate our qualitative observation that the continuous formulation is able to better represent the experimental data by comparing the sum-squared-errors (SSE) between the true parameters used to generate the experimental data to those approximated from formulations \eqref{eq:full_discrete_estimate} and \eqref{eq:full_cont_estimate}. Specifically, we consider the error metrics:
\begin{gather*}
        SSE_{z_{\mu}} := \sum_{i \in \mathcal{I}} (z_{\mu,i} - \bar{z}_{\mu.i})^2 \\
        SSE_{z_{\alpha}} := \sum_{i \in \mathcal{I}} \sum_{j \in \mathcal{I}} (z_{\alpha,ij}-\bar{z}_{\alpha,ij})^2
\end{gather*}
where $\bar{z}_{\mu.i}$ and $\bar{z}_{\alpha,ij}$ denote the actual parameters used in the simulations to generate the experimental data. The results are shown in Table \ref{tab:sse_errors} and demonstrate that the continuous formulation solutions yield significantly smaller sum-squared-errors. Moreover, increased collocation nodes for each finite element are able to reduce the overall error by more than an order of magnitude.

\begin{table}[!htb]
    \centering
    \caption{Sum-squared-errors between the actual and estimated parameters for each formulation and number of collocation nodes per finite element.}
    \begin{tabular}{|c | c | c |} 
        \hline
        Formulation & $SSE_{z_{\alpha}}$ & $SSE_{z_{\mu}}$  \\ [0.5ex]
        \hline \hline
        Discrete, 2 Nodes & $1.02 \times 10^2$ & $4.72 \times 10^{-2}$\\ \hline
        Continuous, 2 Nodes & $2.14 \times 10^1$ & $2.46 \times 10^{-2}$\\ \hline
        Continuous, 4 Nodes & $1.59 \times 10^1$ & $9.05 \times 10^{-4}$\\ \hline
        Continuous, 6 Nodes & $1.59 \times 10^1$ & $9.10 \times 10^{-4}$\\ \hline
    \end{tabular}
    \label{tab:sse_errors}
\end{table}

\section{Conclusions and Future Work} \label{sec:conclusion}

We have presented a unifying  abstraction for representing infinite-dimensional optimization (InfiniteOpt) problems that appear across diverse disciplines. This unifying abstraction introduces the notions of infinite variables (variables parameterized over infinite-dimensional domains). The abstraction also uses measure and differential operators that facilitate the construction of objectives and constraints. The proposed abstraction facilitates knowledge transfer; specifically, it helps identify and transfer modeling constructs across disciplines. For instance, we discussed how chance constraints are analogues of excursion probabilities and how these can be generalized using event constraints; as another example, we show how one can use risk measures to shape time-dependent trajectories in dynamic optimization. The proposed  modeling abstraction aims also to decouple the formulation from transformation schemes (e.g., direct transcription), as we believe that this facilitates analysis and implementation. The proposed abstraction serves as the backbone of a Julia-based modeling framework called \texttt{InfiniteOpt.jl}.
\\

In future work, we will further investigate theoretical crossovers between disciplines. In particular, we are interested in rigorously analyzing the characteristics of risk measures across general infinite domains. Such measures have the potential to enhance the shaping of optimal trajectories in spatio-temporal domains and provide an intuitive analogy for stochastic optimization problems. Event constraints also present an intriguing general constraint class that warrants further research (e.g., connections with disjunctive programming). We also plan to further develop the theoretical and algorithmic foundations for incorporating random field theory into optimization. Furthermore, we will continue enriching the capabilities of \texttt{InfiniteOpt.jl} to support these research endeavors and to make advanced transformation approaches (e.g., MWR approaches) more readily accessible. 

\section*{Acknowledgments}
We acknowledge financial support from the U.S. Department of Energy under grant DE-SC0014114 and from the U.S. National Science Foundation under award 1832208. 

% \appendix
% \section{Some Section}

\bibliography{references}
\end{document}